\documentclass[reqno]{amsart}

\newenvironment{nouppercase}{%
  \renewcommand{\uppercasenonmath}[1]{}}{}

\usepackage{amsmath,amssymb,amsthm,amsfonts,mathrsfs,cite}
\usepackage{array}
\usepackage{fullpage}
\usepackage[colorlinks=true,
linkcolor=blue,
anchorcolor=blue,
citecolor=red
]{hyperref}

\makeatletter
\def\legendre@dash#1#2{\hb@xt@#1{%
  \kern-#2\p@
  \cleaders\hbox{\kern.5\p@
    \vrule\@height.2\p@\@depth.2\p@\@width\p@
    \kern.5\p@}\hfil
  \kern-#2\p@
  }}
\def\@legendre#1#2#3#4#5{\mathopen{}\left(
  \sbox\z@{$\genfrac{}{}{0pt}{#1}{#3#4}{#3#5}$}%
  \dimen@=\wd\z@
  \kern-\p@\vcenter{\box0}\kern-\dimen@\vcenter{\legendre@dash\dimen@{#2}}\kern-\p@
  \right)\mathclose{}}
\newcommand\legendre[2]{\mathchoice
  {\@legendre{0}{1}{}{#1}{#2}}
  {\@legendre{1}{.5}{\vphantom{1}}{#1}{#2}}
  {\@legendre{2}{0}{\vphantom{1}}{#1}{#2}}
  {\@legendre{3}{0}{\vphantom{1}}{#1}{#2}}
}

\newtheorem{theorem}{Theorem}[section]

\newtheorem{corollary}[theorem]{Corollary}

\newtheorem{conjecture}[theorem]{Conjecture}
\theoremstyle{definition}

\theoremstyle{remark}

\newcolumntype{C}[1]{>{\centering\arraybackslash}p{#1}}

\author{Huixi Li}
\address{School of Mathematical Sciences and LPMC, Nankai University, Tianjin 300071, China}
\email{lihuixi@nankai.edu.cn}

\date{\today}

\makeatletter
\@namedef{subjclassname@2020}{\textup{}2020 Mathematics Subject Classification}
\makeatother

\title{Some discussions on the Goldbach conjecture}

\begin{document}
	
\begin{abstract}
According to some discussions based on syllogism, we present results on the binary Goldbach conjecture in three categories: results that are weaker than the Goldbach conjecture, sufficient conditions for the Goldbach conjecture, and results that are similar in nature to the Goldbach conjecture. Additionally, we explore the connections between the Goldbach conjecture and other well-known conjectures.
\end{abstract}
	
\begin{nouppercase}
\maketitle

\hfill \textit{Dedicated to Professor Chen Jingrun on the occasion of the 90th anniversary of his birth}

\hfill \textit{and the 50th anniversary of the birth of Chen's theorem.}

\end{nouppercase}

\tableofcontents

\section{Current status}

In 1742, Goldbach proposed two conjectures: the binary Goldbach conjecture, which states that every even integer greater than $2$ can be expressed as the sum of two primes, and the ternary Goldbach conjecture, which states that every odd integer greater than $5$ can be expressed as the sum of three primes. For a detailed historical account of the conjecture's origins, readers can refer to Section 1 of Vaughan's paper \cite{Vaughan2016}. Note that the ternary Goldbach conjecture can be deduced from the binary Goldbach conjecture, as every odd integer larger than $5$ can be expressed as the sum of $3$ and an even integer greater than $2$. Consequently, the binary Goldbach conjecture is often referred to as the strong Goldbach conjecture, while the ternary Goldbach conjecture is known as the weak Goldbach conjecture.

In 2013, Helfgott \cite{Helfgott2012Minor, Helfgott2013Major, Helfgott2013, Helfgott2015} successfully proved the ternary Goldbach conjecture, providing a conclusive confirmation that every odd integer greater than $5$ can be expressed as the sum of three primes. So our attention will now be focused on the binary Goldbach conjecture for the rest of this paper.

The binary Goldbach conjecture is still unresolved. Nevertheless, there is substantial evidence to suggest that the conjecture is indeed true. Computational verification has been conducted for relatively small even integers, with the binary Goldbach conjecture confirmed for even integers up to $4 \times 10^{18}$ \cite{OHP2014}. While for large even integers, Hardy and Littlewood \cite{HL1923} in 1923 conjectured the following asymptotic formula for the number of representations $r(N)$ of an even integer $N$ as the sum of two primes:
\begin{equation}\label{HLconjecture}
r(N) \sim 2  \prod_{p > 2} \left( 1 - \frac {1} {(p - 1)^2} \right) \times \prod_{\substack{p \mid N \\ p > 2}} \frac {p - 1} {p - 2} \times \frac {N} {(\log N)^2}.
\end{equation}
Notably, the right-hand side of  \eqref{HLconjecture} grows unbounded as $N$ tends to infinity. This implies that large even integers are very likely to have such representations.

While there is substantial evidence suggesting the likely validity of the Goldbach conjecture, it remains an open problem that has yet to be proven or disproven. Understanding the current research status and exploring the vast body of literature surrounding the Goldbach conjecture is of great importance. As of June 2023, a search for the term ``Goldbach" in the MathSciNet database yields over $550$ mathematical publications containing the term, and a more focused search using the Mathematics Subject Classification code ``11P32: Goldbach-type theorems; other additive questions involving primes" as the primary code yields over $800$ publications. Given the extensive research output on the Goldbach conjecture, it becomes meaningful to classify and summarize these articles to provide a valuable resource for researchers interested in this captivating problem. To facilitate a structured approach, we begin the discussion by applying the logical concept of syllogism.

\section{Some discussions based on syllogism}

Suppose we have the following theorem: $A$ implies $B$. Given the additional assumptions that $A^+$ implies $A$ and $A$ implies $A^-$, as well as $B^+$ implies $B$ and $B$ implies $B^-$, we can derive several trivial results:
\begin{itemize}
\item $A$ implies $B$ (as stated in the theorem),
\item $A^+$ implies $B$ (since $A^+$ implies $A$, and $A$ implies $B$),
\item $A$ implies $B^-$ (since $B$ implies $B^-$, and $A$ implies $B$),
\item $A^+$ implies $B^-$ (since $A^+$ implies $A$, and $A$ implies $B^-$).
\end{itemize}
The non-trivial results can be divided into three categories. One may

\noindent 1. improve the theorem:
\begin{itemize}
    \item $A^-$ implies $B$,
    \item $A$\,\,\,\, implies $B^+$,
    \item $A^-$ implies $B^+$,
\end{itemize}
2. obtain a stronger result under a stronger assumption:
\begin{itemize}
    \item $A^+$ implies $B^+$,
\end{itemize}
3. or obtain a weaker result under a weaker assumption:
\begin{itemize}
    \item $A^-$ implies $B^-$.
\end{itemize}
In addition, it is possible to obtain results of a similar nature to the original theorem:
\begin{itemize}
    \item $A'$ \, implies $B'$.
\end{itemize}
Suppose we have an unsolved conjecture: $C$ implies $D$. Given the additional assumptions that $C^+$ implies $C$ and $C$ implies $C^-$, as well as $D^+$ implies $D$ and $D$ implies $D^-$, since an unsolved conjecture is not proved or disproved, the following 18 types of results are all non-trivial:
\begin{itemize}
    \item \color{red}$C$ \,\, implies $D$, \color{black}
    \item \color{red}$C^-$ implies $D$, \color{black}
    \item \color{yellow}$C^+$ implies $D$, \color{black}
    \item \color{red}$C$ \,\, implies $D^+$, \color{black}
    \item \color{red}$C^-$ implies $D^+$, \color{black}
    \item \color{yellow}$C^+$ implies $D^+$, \color{black}
    \item \color{blue}$C$ \,\, implies $D^-$, \color{black}
    \item \color{blue}$C^-$ implies $D^-$, \color{black}
    \item \color{blue}$C^+$ implies $D^-$, \color{black}
    \color{red}
    \item $C$ \,\, may not imply $D$,
    \item $C^-$ may not imply $D$,
    \item $C^+$ may not imply $D$,
    \item $C$ \,\, may not imply $D^+$,
    \item $C^-$ may not imply $D^+$,
    \item $C^+$ may not imply $D^+$,
    \item $C$ \,\, may not imply $D^-$,
    \item $C^-$ may not imply $D^-$,
    \item $C^+$ may not imply $D^-$.
    \color{black}
\end{itemize}
Regarding the Goldbach conjecture, it remains unproven, and no counterexamples have been found. Therefore, the types of results marked in red above are not applicable to the Goldbach conjecture.  However, there have been notable advancements marked in blue, indicating partial progress towards the conjecture. Furthermore, there are results marked in yellow that establish sufficient conditions for the Goldbach conjecture.

Finally, one may obtain results that share a similar nature with the conjecture:

\begin{itemize}
    \item $C'$ \, implies $D'$.
\end{itemize}

Based on these discussions, we can classify the current results on the binary Goldbach conjecture into three categories: results that are weaker than the Goldbach conjecture, sufficient conditions for the Goldbach conjecture, and results that share a similar nature with the Goldbach conjecture. These categories will be discussed in detail in Sections \ref{Some results weaker than the Goldbach conjecture}, \ref{some sufficient conditions for the Goldbach conjecture}, and \ref{some results of the same flavor as the Goldbach conjecture}, respectively. Additionally, in Section \ref{connectionsec}, we will explore the connections between the Goldbach conjecture and other renowned conjectures, such as the twin prime conjecture and the Riemann hypothesis. By organizing the current research findings in this manner, we aim to provide a comprehensive overview of the progress made in understanding the binary Goldbach conjecture.

The author acknowledges that the presented overview of the Goldbach conjecture may not encompass all relevant contributions made by researchers in the field. We apologize for any unintentional omissions and recognize the significant contributions made by numerous researchers to the study of the Goldbach conjecture. The intent is to provide a general summary and understanding of the topic while acknowledging the vastness of the research landscape surrounding the Goldbach conjecture.

\section{Some results weaker than the Goldbach conjecture}\label{Some results weaker than the Goldbach conjecture}

In this section, we will discuss several results that are weaker than the Goldbach conjecture but still provide valuable insights into the problem. These results include the representation of positive even integers as the sum of a prime and almost primes, the representation of positive integers as the sum of a prime and a square-free number, the exceptional set of Goldbach numbers, the distributions of Goldbach numbers and exceptional Goldbach numbers, the ternary Goldbach conjecture with restrictions, and the Linnik-Goldbach representation.

\subsection{Goldbach conjecture and almost primes}

By the sieve method we can naturally obtain results on almost primes. A natural number is called a $k$-almost prime or an almost prime of order $k$ if it has $k$ prime divisors counted with multiplicity. Clearly primes and $1$-almost primes are the same, so the Goldbach conjecture is also known as $1 + 1$, while the $a+b$ results are partial results towards the Goldbach conjecture, which means every sufficiently large even integer can be written as the sum of two almost primes, each with at most $a$ and $b$ prime divisors respectively counted with multiplicity. Another way to express this is $N = P_a + P_b$ is solvable for all sufficiently large even integers, where $P_r$ means an almost primes with at most $r$ prime divisors counted with multiplicity. Since Brun \cite{Brun1920} first proved $9+9$ in 1920, mathematicians have been working on representing large even integers as the sum of two almost primes over the last 100 years. Currently the best result is due to Chen \cite{Chen1973}, who proved $1 + 2$ in 1973. For a collection of related papers we refer the readers to Wang's book \cite{Wang1982}.

\begin{table}[ht]\label{a+b table}
\begin{center}
\begin{tabular}{ |C{6cm}|C{6cm}|C{1cm}| }
\hline
\text{Brun}
& $9+9$
&
1920 \\
\hline
\text{Rademacher}
& $7+7$
&
1924 \\
\hline
\text{Estermann}
& $6+6$
&
1932 \\
\hline
\text{Buchstab}
& $5+5$
&
1938 \\
\hline
\text{Buchstab}
& $4+4$
&
1940 \\
\hline
\text{Renyi}
& $1+c$
&
1948 \\
\hline
\text{Vinogradov}
& $3 + 3$
&
1956 \\
\hline
\text{Wang}
& $2 + 3$
&
1957 \\
\hline
\text{Barban and Pan}
& $1+5$
&
1962 \\
\hline
\text{Barban, Pan, and Wang}
& $1+4$
&
1962 \\
\hline
\text{Buchstab, Bombieri, and Vinogradov}
& $1+3$
&
1965 \\
\hline
\text{Chen}
& $1+2$
&
1973 \\
\hline
\end{tabular}
\end{center}
\caption{Sieve theory results of type $a + b$}
\label{table1}
\end{table}
More explicitly, Chen proved the following result in 1973.
\begin{theorem}[Chen, 1973]
Let $P_x(1, 2)$ be the number of ways to write a sufficiently large even integer $x$ as the sum of a prime and an almost prime with at most $2$ prime divisors, then we have
\begin{align}\label{chenconstanteq}
P_x(1, 2) \geq \frac {0.67 x C_x} {(\log x)^2},
\end{align}
where
\begin{align*}
C_x =
\prod_{p > 2} \left( 1 - \frac {1} {(p - 1)^2} \right)
\times \prod_{\substack{p \mid x \\ p > 2}} \frac {p - 1} {p - 2}.
\end{align*}
Therefore, every sufficiently large even integer can be written as the sum of a prime and product of at most two primes.
\end{theorem}
Then constant term $0.67$ in~\eqref{chenconstanteq} improves his earlier result \cite{Chen1966} announced in 1966, in which Chen obtained a smaller constant $0.098$ in the lower bound estimate. There are further improvements on the constant term in Chen's theorem. Currently the record belongs to Wu \cite{Wu2008}, who proved the following theorem in 2008, improving previous results of Chen \cite{Chen1966, Chen1973}, Halberstam \cite{Halberstam1975}, Chen \cite{Chen1978_1}, Chen \cite{Chen1978_2}, Cai-Lu \cite{CL2002}, Wu \cite{Wu2004}, and Cai \cite{Cai2008}.

\begin{theorem}[Wu, 2008]
We have $P_x(1, 2) \geq \frac{0.899 x C_x}{(\log x)^2}$.
\end{theorem}

\begin{table}[ht]\label{Constant term in Chen's theorem table}
\begin{center}
\begin{tabular}{ |C{6cm}|C{6cm}|C{1cm}| }
\hline
\text{Chen}
& $P_x(1, 2) \geq \frac{0.67 x C_x}{(\log x)^2}$
&
1973 \\
\hline
\text{Halberstam}
& $P_x(1, 2) \geq \frac{0.689 x C_x}{(\log x)^2}$
&
1975 \\
\hline
\text{Chen}
& $P_x(1, 2) \geq \frac{0.7544 x C_x}{(\log x)^2}$
&
1978 \\
\hline
\text{Chen}
& $P_x(1, 2) \geq \frac{0.81 x C_x}{(\log x)^2}$
&
1978 \\
\hline
\text{Cai-Lu}
& $P_x(1, 2) \geq \frac{0.8285 x C_x}{(\log x)^2}$
&
2002 \\
\hline
\text{Wu}
& $P_x(1, 2) \geq \frac{0.836 x C_x}{(\log x)^2}$
&
2004 \\
\hline
\text{Cai}
& $P_x(1, 2) \geq \frac{0.867 x C_x}{(\log x)^2}$
&
2008 \\
\hline
\text{Wu}
& $P_x(1, 2) \geq \frac{0.899 x C_x}{(\log x)^2}$
&
2008 \\
\hline
\end{tabular}
\caption{Constant term in Chen's theorem}
\label{table2}
\end{center}
\end{table}

Note that Chen's theorem holds for \emph{sufficiently large} even integers. What is the specific threshold for the validity of Chen's representation for positive even integers? In 2022, Bordignon, Johnston, and Starichkova \cite{BJS2022, Bordignon2022} obtained a result that improves upon Yamada's earlier work \cite{Yamada2015} in 2015. Bordignon, Johnston, and Starichkova's proof incorporates \emph{Siegel zeros}, a concept defined in \cite[Chapter 21]{Davenport2000}, which prompts the intriguing question of whether Zhang's new results on Siegel zeros \cite{Zhang2022} can further enhance the estimates in the explicit version of Chen's theorem.

\begin{theorem}[Bordignon-Johnston-Starichkova, 2022]
Every even number greater than $\exp(\exp(32.6))$ can be represented as the sum of a prime and a product of at most two primes.
\end{theorem}
In 2001, Richstein \cite{Richstein2001} found that the smaller primes in the representations of even integers up to $n = 10^{14}$ as the sum of two primes do not exceed $5569$. In 1988, Granville, van de Lune, and te Riele \cite{GVT1988} conjectured the order of the smaller primes in the Goldbach representations for even integers up to $n$ is $O((\log n)^2 \log \log n)$. In the setting of Chen's theorem with small primes, Cai proved the following theorem, which surpasses the previous results of Cai \cite{Cai2002} and Li-Cai \cite{CL2011}.
\begin{theorem}[Cai, 2014]
Let $N$ be a sufficiently large even integer, then $N = p + P_2$ is solvable for some $p \leq N^{0.941}$.
\end{theorem}
\begin{table}[ht]\label{Chen theorem with small primes table}
\begin{center}
\begin{tabular}{ |C{6cm}|C{6cm}|C{1cm}| }
\hline
\text{Cai}
& $p \leq N^{0.95}$
&
2002 \\
\hline
\text{Li-Cai}
& $p \leq N^{0.945}$
&
2011 \\
\hline
\text{Cai}
& $p \leq N^{0.941}$
&
2015 \\
\hline
\end{tabular}
\end{center}
\caption{Chen's theorem with small primes}
\label{table3}
\end{table}

Cai obtained a similar result \cite{Cai2008} in 2008, in which it is required that both $p$ and $P_2$ are close to $N/2$. This improves the previous results of Ross \cite{Ross1978}, Wu \cite{Wu1993}, Wu \cite{Wu1994}, Salerno-Vitolo \cite{SV1993}, Cai-Lu \cite{CL1999}, and Wu \cite{Wu2004}.
\begin{theorem}[Cai, 2008]
Let $N$ be a sufficiently large even integer and let $U = N^{0.97}$. Then $N = p + P_2$ with $N/2 - U \leq p, P_2 \leq N/2 + U$ is solvable.
\end{theorem}
\begin{table}[ht]
\begin{center}
\begin{tabular}{ |C{6cm}|C{6cm}|C{1cm}| }
\hline
\text{Ross}
& $U = N^{0.98}$
&
1978 \\
\hline
\text{Wu}
& $U = N^{0.974}$
&
1993 \\
\hline
\text{Wu}
& $U = N^{0.973}$
&
1994 \\
\hline
\text{Salerno-Vitolo}
& $U = N^{0.9729}$
&
1993 \\
\hline
\text{Cai-Lu}
& $U = N^{0.972}$
&
1999 \\
\hline
\text{Wu}
& $U = N^{0.971}$
&
2004 \\
\hline
\text{Cai}
& $U = N^{0.97}$
&
2008 \\
\hline
\end{tabular}
\end{center}
\caption{Chen's theorem with almost equal summands}
\label{table4}
\end{table}

Here is another version of Chen's theorem with restrictions by Cai \cite{Cai2017} in 2017, improving the work of Heath-Brown and Li \cite{HBL2016} in 2016.
\begin{theorem}[Cai, 2017]
Let $N$ be a sufficiently large even integer. Then $N = p + P_2$ with $p + 6 = P_{14}$ is solvable.
\end{theorem}

Finally, by utilizing the upper bound sieve, Bombieri and Davenport \cite{BD1966} proved in 1966 that the number of solutions of the equation $N = p_1 + p_2$ with prime variables $p_1$ and $p_2$ for sufficiently large even integers is at most $(8 + \epsilon) \frac{N C_N}{(\log N)^2}$, where $\epsilon$ is a small positive constant, while the correct order for the number of such representations was obtained earlier by Brun \cite{Brun1920} in 1920. Chen \cite{Chen1978_3} subsequently reduced the constant $8 + \epsilon$ to $7.8342$ in 1978, and Wu \cite{Wu2004} further improved Chen's result in 2004, achieving a constant of $7.8209$. More generally, Kan studied the upper bound and lower bound for number of solutions of $N = P_s + P_r$ for sufficiently large even integers $N$ and positive integers $s \geq 1$ and $r \geq 2$ in the 1990s \cite{Kan1991, Kan1992}, in which he obtained a lower bound whose order misses the presumably correct upper bound order by a factor of $\log \log N$.

\subsection{Representing positive integers as the sum of a prime and a square-free number}
Square-free numbers are positive integers not divisible by prime squares. The set of square-free numbers includes $1$, $2$, $3$, $5$, $6$, $7$, $10$, $11$, $13$, $\cdots\cdots$. Unlike the primes and the almost primes, the square-free numbers have a positive density of $\frac{6}{\pi^2}$ among the positive integers. In 1931 Estermann proved the following theorem, which can be viewed as a partial result towards the Goldbach conjecture \cite{Estermann1931}, as primes are clearly square-free.
\begin{theorem}[Estermann, 1931]
Every sufficiently large integer may be written as the sum of a prime and a square-free number.
\end{theorem}

Recently, there has been considerable attention given to the representation of positive integers as the sum of a prime and a square-free number, which we refer to Estermann's representation, especially after Dudek's proof \cite{Dudek2017} of the explicit version of Estermann's theorem in 2017.

\begin{theorem}[Dudek, 2017]
Every integer greater than $2$ may be written as the sum of a prime and a square-free number.
\end{theorem}

Researchers have also studied Estermann's representation with restrictions. For example, Yau \cite{Yau2021} obtained an estimate in 2021 for the weighted number of ways a sufficiently large integer can be represented as the sum of a prime congruent to $a$ modulo $q$ and a square-free integer. Similarly, Hathi and Johnston \cite{HJ2021} investigated Estermann's theorem with divisibility conditions on the square-free number part in 2021, building upon previous work by Francis and Lee \cite{FrancisLee2020, FrancisLee2022}.

Similar to Cai's requirement for smaller primes in Chen's representation, researchers have explored the possibility of imposing a restriction that the primes in Estermann's representation can be taken to be small. Notably, Dalton and Trifonov \cite{DT2023} made significant progress in this direction and proved the following theorem in 2023.

\begin{theorem}[Dalton-Trifonov,  2023]
Every positive integer $n$ which is not equal to $1$, $2$, $3$, $6$, $11$, $30$, $155$, or $247$ can be written as the sum of a prime and a square-free number, in which the prime does not exceed $\sqrt{n}$.
\end{theorem}

A natural question arises: Can we express a sufficiently large integer as the sum of a prime and a square-free number with a small number of prime divisors? This line of research combines the representations of Estermann and Chen, making it an intriguing problem to explore.

If one is familiar with the proof of Chen's theorem that uses the sieve method, it is known that Chen's $1+2$ result can actually be strengthened in its statement. More explicitly, in order to sieve out numbers with many prime factors, Chen \cite{Chen1973} considers numbers of the form $N-p$ whose prime factors are greater than $N^{1/10}$. This ensures that the count of such numbers that are not squaere-free does not exceed $N^{9/10+\epsilon}$ for any small positive constant $\epsilon$. This point is made more clearly in Ross's article \cite{Ross1975}. Hence, Chen has effectively resolved the aforementioned problem for sufficiently large even numbers.

Now, does a similar conclusion hold for sufficiently large odd numbers? The answer is affirmative. Utilizing the theorem of Pan, Ding, and Wang \cite{PDW1975}, we know that every sufficiently large odd number $M$ can be expressed as the sum of a prime and either twice an odd prime or twice the product of two odd primes. Similarly, we can impose the additional condition that the odd primes in the Pan-Ding-Wang representation are greater than or equal to $M^{1/10}$, since the count of such numbers is much smaller than the main term, it is permissible to further restrict the second term in the Pan-Ding-Wang representation to be square-free. More explicitly, Li  \cite{Li2019} proved in 2019 that:
\begin{theorem}[Pan-Ding-Wang, 1975, L., 2019]\label{Li2019thm}
Any sufficiently large odd number can be written as the sum of a prime and $2$ times an almost prime with at most two prime divisors. Let $R(M)$ be the number of primes $p$ up to $M$ such that $M - p$ is $2$ times and odd prime or $2$ times the product of two distinct odd primes, then for sufficiently large odd integers $M$ we have
\begin{align*}
R(M) \geq \frac {0.32 M C_M} {(\log M)^2}.
\end{align*}
\end{theorem}
By combining this result with Chen's theorem, we have successfully addressed the problem posed above in its entirety.
\begin{corollary}
Every sufficiently large integer can be written as the sum of a prime and a square-free number with at most three prime divisors.
\end{corollary}

Can we further improve the constraint on the number of prime factors for the square-free number in the aforementioned conjecture to at most two? This question is related to the Lemoine conjecture proposed in 1895.
\begin{conjecture}[Lemoine, 1895]
Any odd integer greater than $5$ can be expressed as the sum of a prime and twice another prime.
\end{conjecture}
Therefore, reducing the constraint from three to two would require proving the validity of the Lemoine conjecture for sufficiently large odd numbers. This is an extremely challenging problem and can be considered even more difficult than the Goldbach conjecture.

It is worth noting that Chen's theorem provides a lower bound estimate for the number of solutions to the Diophantine equation $N = p + P_2$ for sufficiently large even integers $N$, where $p$ is a prime variable and $P_2$ is a variable that is either a prime or the product of two primes. In a similar vein, Theorem~\ref{Li2019thm} gives a lower bound estimate for the number of solutions to the Diophantine equation $M = p + 2P_2$ for sufficiently large odd integers $M$. This prompts us to ask whether we can provide a lower bound estimate for the number of solutions to the Diophantine equation $M = 2p + P_2$ for sufficiently large odd integers $M$. Furthermore, can we generalize this to the equation $N = a p + b P_2$ for given positive integers $a$ and $b$, subject to certain conditions on $N$? The answer is affirmative. In 2023, Li \cite{Li2023} established the following result, providing a lower bound estimate for the number of solutions to such Diophantine equations.
\begin{theorem}[L., 2023]\label{Li2023thm}
For two relatively prime square-free positive integers $a$ and $b$, let $N$ be a sufficiently large integer that is relatively prime to both $a$ and $b$, and let $N$ be even if $a$ and $b$ are both odd. Let $R_{a, b}(N)$ be the number of primes $p$ such that $a p$ and $N - a p$ are both square-free, $b \mid (N - a p)$, and $\frac{N - a p}{b}$ is either a prime or the product of two primes. We have
\begin{align*}
R_{a, b}(N)
\geq
\frac{0.68 C_N N}{a b (\log N)^2}
\times
\prod_{\substack{p \mid a b \\ p > 2}}
\frac{p - 1}{p - 2}.
\end{align*}
\end{theorem}
Indeed, Theorem~\ref{Li2023thm} can be extended to a more general setting by relaxing the condition that $a$ and $b$ are square-free. In this case, we define $R_{a, b}(N)$ as the number of primes $p$ for which $\frac{N - ap}{b}$ is either a prime or the product of two distinct primes coprime to $b$. By carefully modifying the details of the proof, we can obtain a lower bound estimate for the number of representations of a given integer satisfying certain conditions as the sum of $a p$ and $b P_2$, where $p$ and $P_2$ are selected from specific arithmetic progressions. Finally, we remark that Li \cite{Li2023} also proved the number of solutions of the Diophantine equation $N = a p_1 + b p_2$ with prime variables $p_1$ and $p_2$ has a satisfactory upper bound.

\subsection{Exceptional set of Goldbach numbers}\label{exceptionalsetofGoldbachnumberssec}

The circle method utilizes complex analysis and number-theoretic techniques to investigate problems related to additive number theory. In the context of the Goldbach conjecture, the circle method allows us to analyze the exceptional set of Goldbach numbers, which consists of positive even integers that cannot be expressed as the sum of two primes. By applying the circle method, we can obtain estimations for the cardinality of this exceptional set.

Let $E(N)$ be the number of positive even integers up to $N$ that are not the sum of two primes. Then the Goldbach conjecture is equivalent to that for all $N \geq 2$, we have $E(N) = 1$. While we do not have a proof for the statement $E(N) = 1$ for all $N \geq 2$, we can establish weaker results. Assuming the generalized Riemann hypothesis (GRH), Hardy and Littlewood \cite{HL1924} proved $E(N) \ll N^{1/2 + \epsilon}$ in 1924. Unconditionally, the best known bound for $E(N)$ is due to Pintz \cite{Pintz2018_1, Pintz2018_2, Pintz2023}, who proved that $E(N) \ll N^{0.72}$ in 2018, improving previous bounds obtained by Chudakov \cite{Chudakov1937}, van der Corput \cite{vanderCorput1937}, Estermann \cite{Estermann1938}, Vaughan \cite{Vaughan1972}, Montgomery-Vaughan \cite{MV1975}, Chen-Pan \cite{CP1980}, Chen \cite{Chen1983}, Chen-Liu \cite{CL1989}, Li \cite{Li1999}, Li \cite{Li2000exp}, and Lu \cite{Lu2010}.

\begin{theorem}[Pintz, 2018]
Let $E(N)$ be the number of positive even integers up to $N$ that are not the sum of two primes. We have $E(N) \ll N^{0.72}$.
\end{theorem}

\begin{table}[ht]
\begin{center}
\begin{tabular}{ |C{6cm}|C{6cm}|C{1cm}| }
\hline
\text{Chudakov}
& $E(N) \ll \frac{N}{(\log N)^A}$
&
1937 \\
\hline
\text{van der Corput}
& $E(N) \ll \frac{N}{(\log N)^A}$
&
1937 \\
\hline
\text{Estermann}
& $E(N) \ll \frac{N}{(\log N)^A}$
&
1938 \\
\hline
\text{Vaughan}
& $E(N) \ll \frac{N}{e^{c \sqrt{\log N}}}$
&
1972 \\
\hline
\text{Montgomery-Vaughan}
& $E(N) \ll N^{1 - \delta}$
&
1975 \\
\hline
\text{Chen-Pan}
& $E(N) \ll N^{0.99}$
&
1980 \\
\hline
\text{Chen}
& $E(N) \ll N^{0.96}$
&
1983 \\
\hline
\text{Chen-Liu}
& $E(N) \ll N^{0.95}$
&
1989 \\
\hline
\text{Li}
& $E(N) \ll N^{0.921}$
&
1999 \\
\hline
\text{Li}
& $E(N) \ll N^{0.914}$
&
2000 \\
\hline
\text{Lu}
& $E(N) \ll N^{0.879}$
&
2010 \\
\hline
\text{Pintz}
& $E(N) \ll N^{0.72}$
&
2018 \\
\hline
\end{tabular}
\end{center}
\caption{Exceptional set of Goldbach numbers}
\label{table5}
\end{table}

We remark that Montgomery and Vaughan \cite{MV1975} made significant contributions to the study of the exceptional Goldbach numbers in their influential work published in 1975, in which they investigated the two cases whether the Siegel zeros exist or do not exist and obtained power saving in the estimate of $E(N)$.

Additionally, it is worth mentioning that the exceptional set in Goldbach's problem in short intervals has been studied by many mathematicians. Grimmelt \cite{Grimmelt2022} proved in 2022 that $E(X, H) - E(X) \ll H (\log X)^{-A}$ in the range $H > X^{11/180 + \epsilon}$ for any $A > 0$ and $\epsilon > 0$, improving previous results by Ramachandra \cite{Ramachandra1973}, Mikawa \cite{Mikawa1992}, Perelli-Pintz \cite{PP1993}, Dufner \cite{Dufner1994, Dufner1995}, Jia \cite{Jia1995short, Jia1995short2, Jia1996short}, Li \cite{Li1995}, and Harman \cite{Harman2007}. With respect to power saving, Matom\"{a}ki \cite{Matomaki2008} proved in 2008 that for $H = X^\theta$ with $1/5 < \theta \leq 1$, there exists some $\delta_0 = \delta_0(\theta) > 0$ such that $E(X, H) - E(X) \ll H^{1 - \delta_0}$, improving previous work by Luo-Yao \cite{LY1981}, Yao \cite{Yao1982}, Peneva \cite{Peneva2001}, and Languasco \cite{Languasco2004Goldbach}. Also, see \cite[Theorem 1.3]{MRT2019} in 2019 by Matom\"{a}ki, Radziwi\l \l, and Tao for the averaged Goldbach representations in short intervals.

Furthermore, we remark that Li \cite{Li2003execptionalset} in 2003 considered the linear equation $N = a_1 p_1 + a_2 p_2$ with coprime natural numbers coefficients $a_1$ and $a_2$ and with prime variables $p_1$ and $p_2$. Let $E(N)$ be the number of integers $m$ with $(1 - \epsilon) N \leq m \leq N$, where $\epsilon$ is a sufficiently small positive number, and $m$ is coprime to $a_1 a_2$, with the additional condition that $m$ is even if $a_1 a_2$ is odd. Then Li proved the estimate $E(N) \ll N^{1 - \delta}$ for some $\delta > 0$ for all $N \gg \max(a_1, a_2)^{31}$, providing an explicit version of Liu's result \cite[Theorem 2]{Liu1987} in 1987.

Finally, the exceptional set corresponding to the binary Goldbach problem $N = p_1 + p_2$ with prime summands $p_1$ and $p_2$ in arithmetic progressions with large moduli has also been studied by Liu-Zhan \cite{LZ1997}, Bauer \cite{Bauer2012, Bauer2017, Bauer2017_2}, Bauer-Wang \cite{BW2013}, Martin \cite{Martin2022}, and Salmenssu \cite{Salmensuu2022}. See Salmenssu's paper \cite{Salmensuu2022} for more details. Chen's theorem $N = p + P_2$ with summands $p$ and $P_2$ in arithmetic progressions with large moduli was studied by Cai and Lu \cite{CL1999_2} in 1999.

\subsection{Distributions of Goldbach numbers and exceptional Goldbach numbers}
We have discussed the cardinality of the set of exceptional Goldbach  numbers in the previous section. It is known that the Goldbach numbers have full density among the positive even integers, while the exceptional Goldbach numbers have density zero. Even though we have only found just one exceptional Goldbach numbers so far, an interesting question is: how are they distributed? If there are many exceptional Goldbach numbers, do they tend to cluster together or are they evenly distributed? One way to look at this is to compute the moment of the gaps between consecutive Goldbach numbers.

Let $g_n$ be the $n$-th Goldbach number. If $X \geq 4$ is a positive integer, then the Goldbach conjecture is equivalent to
\[
\sum_{g_n \leq x}
(g_{n + 1} - g_n)^\gamma
= 2^{\gamma} \left\lfloor \frac{X - 1}{2} \right\rfloor = 2^{\gamma - 1} X + O(1)
\]
for any positive constant $\gamma$. It is hard to prove this strong result, however, we can prove weaker results in the sense that one may relax the error term estimation $O(1)$ while putting restrictions on the size of $\gamma$. Pintz \cite{Pintz2018_1, Pintz2018_2, Pintz2023} proved the following result in 2018, improving a previous result of Mikawa \cite{Mikawa1993} in 1993.
\begin{theorem}[Pintz, 2018, Theorem C]
Let $g_n$ be the $n$-th Goldbach number. For $\gamma < \frac{341}{21}$ we have
\begin{equation}\label{moment sum of gaps}
\sum_{g_n \leq X}
(g_{n + 1} - g_n)^\gamma
= 2^{\gamma - 1} X
+ O(X^{1 - \delta}),
\end{equation}
where $\delta > 0$ is a small positive constant.
\end{theorem}

\begin{table}[ht]
\begin{center}
\begin{tabular}{ |C{6cm}|C{6cm}|C{1cm}| }
\hline
\text{Mikawa}
&
$\gamma < 3$
&
1993 \\
\hline
\text{Pintz} &
$\gamma < \frac{341}{21}$
&
2018 \\
\hline
\end{tabular}
\end{center}
\caption{Moments of gaps between Goldbach number}
\label{table6}
\end{table}

One way to understand Equation \eqref{moment sum of gaps} is: when the error term is fixed to be $O(X^{1 - \delta})$, then the larger values of $\gamma$ can be taken, the stronger the result becomes, which means the gaps between consecutive Goldbach numbers can not be too large, i.e., the exceptional Goldbach numbers tend not to be clustered.

Pintz \cite{Pintz2018_1, Pintz2018_2, Pintz2023} has also proved that most small intervals do not contain a large number of exceptional Goldbach numbers. Specifically, it has been shown that the number of small intervals of the form $(n, (n + \log n)^2]$ that contains at most a constant number of exceptional Goldbach numbers has density zero.

\begin{theorem}[Pintz, 2018, Theorem E]
There are explicitly calculable absolute constants $K$ and $C$ such that for all but $C X^{3/5} (\log X)^{12}$ numbers $n \leq X$  we have
\[
E(n + (\log n)^2) - E(n) \leq K.
\]
\end{theorem}

Harman \cite{Harman2020} considered the Diophantine approximation problem for the Goldbach numbers and proved the following result in 2020.
\begin{theorem}[Harman, 2020]
Given $\epsilon > 0$, for every irrational $\alpha$ and arbitrary real $\beta$ there are infinitely many solutions to
\[
|| \alpha n + \beta || < n^{- \frac{5}{6}},
\]
where $n$ is required to be a Goldbach number and $||x||$ means the distance form $x$ to a nearest integer.
\end{theorem}
If the Goldbach conjecture is true, then the above inequality holds for infinitely many Goldbach numbers $n$ when we replace $n^{-\frac{5}{6}}$ by $2 n^{-1}$.

\subsection{Ternary Goldbach conjecture with restrictions}
The ternary Goldbach conjecture has been proved, i.e., every odd integer greater than $5$ is the sum of three primes. If we can further restrict that the smallest prime equals $3$ in the representation of an odd integer as the sum of three primes, then this would imply the binary Goldbach conjecture is also true. There are results weaker than this statement, in which the smallest prime in the above representation may not be as small as $3$. Under the assumption of the GRH holds, Montgomery and Vaughan \cite{MV1975} proved the following result in 1975, improving Linnik's result \cite{Linnik1952} in 1952.
\begin{theorem}[Montgomery-Vaughan, 1975]
Suppose the GRH is true. Let $N$ be a sufficiently large odd integer. The equation
$N = p_1 + p_2 + p_3$ is solvable in primes $p_1$, $p_2$, and $p_3$, where $p_1 \ll (\log X)^2$.
\end{theorem}
Unconditionally, Cai \cite{Cai2013} proved the following result in 2013, improving previous results of Pan \cite{Pan1959}, Montgomery-Vaughan \cite{MV1975}, Zhan \cite{Zhan1995}, Wong \cite{Wong1996}, and Jia \cite{Jia1996}.

\begin{theorem}[Cai, 2013]
Let $N$ be a sufficiently large odd integer. The equation
$N = p_1 + p_2 + p_3$ is solvable in primes $p_1$, $p_2$, and $p_3$, where $p_1 \leq N^{\frac{11}{400} + \epsilon}$.
\end{theorem}

\begin{table}[ht]
\begin{center}
\begin{tabular}{ |C{6cm}|C{6cm}|C{1cm}| }
\hline
\text{Pan}
& $p_1 \leq  N^{\frac{1}{4} + \epsilon}$
&
1959 \\
\hline
\text{Montgomery-Vaughan}
& $p_1 \leq  N^{\frac{7}{72} + \epsilon}$
&
1975 \\
\hline
\text{Zhan}
& $p_1 \leq  N^{\frac{7}{120} + \epsilon}$
&
1995 \\
\hline
\text{Wong}
& $p_1 \leq  N^{\frac{7}{216} + \epsilon}$
&
1996 \\
\hline
\text{Jia}
& $p_1 \leq  N^{\frac{7}{240} + \epsilon}$
&
1996 \\
\hline
\text{Cai}
& $p_1 \leq  N^{\frac{11}{400} + \epsilon}$
&
2013 \\
\hline
\end{tabular}
\end{center}
\caption{Ternary Goldbach conjecture with small primes}
\label{table7}
\end{table}

As discussed in Section~\ref{exceptionalsetofGoldbachnumberssec}, similar to the estimation of the exceptional set of Goldbach numbers, we can also estimate the size of the set of ``exceptional $C$-odd numbers": for a given constant $C$, an odd number is considered an \emph{exceptional $C$-odd number} if it cannot be expressed as the sum of three primes, with the smallest prime being smaller than $C$.
More explicitly, Pintz \cite{Pintz2018_1, Pintz2018_2, Pintz2023} proved the following result in 2018.
\begin{theorem}[Pintz, 2018, Theorem B]
All but $O(N^{3/5} (\log N)^{10})$ odd numbers can be written as the sum of three primes with one prime less than $C$, a given absolute constant.
\end{theorem}

In addition to the aforementioned results, there are many other results that impose restrictions on the ternary Goldbach conjecture. These include the almost equal variables version \cite{Haselgrove1951, Pan1959, Chen1965, PP1989, PP1990, Jia1989, Jia1991_1, Jia1991_2, Jia1991_3, Jia1992, Jia1991_4, Jia1994, Zhan1991, BH1998, MMS2017, LS2013}, almost primes version \cite{Tolev1999, Tolev2000_1, Tolev2000_2, Peneva2000, Meng2007, Meng2009, GT2022}, density version \cite{Shen2016, Shao2014, LP2010}, Fouvry-Iwaniec primes version \cite{Grimmelt2022FIprimes}, primes in the Beatty sequences $\{\lfloor \alpha n + \beta \rfloor\}_{n = 1}^\infty$ version \cite{BGN2007, Kumchev2008, LS2013, Vaughan2014, SC2021}, primes in the Piatetski-Shapiro sequences $\{\lfloor n^c \rfloor\}_{n = 1}^\infty$ version \cite{BF1992, Rivat1992, Jia1995, Kumchev1997, Liu1998, LZ2022}, primes with given primitive roots version \cite{FKS2021}, small prime solutions with coefficients version \cite{Baker1967, LT1989, LT1991, Choi1997, LW1998, Li2001smallprime, LT2005, CK2006, CT2017}, and many other versions.

\subsection{Linnik-Goldbach representation}\label{LinnikGoldbach}

The set of integers that can be expressed as the sum of at most $K$ powers of $2$ is very sparse. Consequently, an alternative method to approach the Goldbach conjecture is to express large even numbers as the sum of two primes, along with at most $K$ powers of $2$ for some constant $K$. If the value of $K$ can be taken to be $0$, meaning that any sufficiently large even number can be expressed as the sum of two primes without any powers of $2$, then the Goldbach conjecture would hold for sufficiently large even numbers. Indeed, we cannot prove such a strong statement. However, Linnik \cite{Linnik1953} proved in 1953 that there exists a large enough constant $K$ such that every sufficiently large even integer can be written as the sum of two primes and at most $K$ powers of $2$. In 1975, Gallagher \cite{Gallagher1975} used a different method to prove the same result. In 1998, Liu-Liu-Wang \cite{LLW1998_2} proved that the value of $K$ can be specifically taken as $54000$. Since then, mathematicians have made a series of improvements to this constant. These include the works of Li \cite{Li2000Linnik}, Wang \cite{Wang1999}, Li \cite{Li2001}, Heath-Brown-Puchta \cite{HBP2002}, Pintz-Puzsa \cite{PR2003}, Liu-Lv \cite{LL2011}, and Pintz-Ruzsa \cite{PR2020}. In 2020, Pintz and Ruzsa \cite{PR2020} proved the following theorem.

\begin{theorem}[Pintz-Ruzsa, 2020]
Every sufficiently large even integer is the sum of two primes and at most $8$ powers of $2$.
\end{theorem}

\begin{table}[ht]
\begin{center}
\begin{tabular}{ |C{6cm}|C{6cm}|C{1cm}| }
\hline
\text{Linnik}
&
$K$
&
1953 \\
\hline
\text{Gallagher}
&
$K$
&
1975 \\
\hline
\text{Liu-Liu-Wang}
&
54000
&
1998 \\
\hline
\text{Li}
&
25000
&
2000 \\
\hline
\text{Wang}
&
2250
&
1999 \\
\hline
\text{Li}
&
1906
&
2001 \\
\hline
\text{Heath-Brown--Puchta}
&
13
&
2002 \\
\hline
\text{Pintz-Ruzsa}
&
13
&
2003 \\
\hline
\text{Liu-Lv}
&
12
&
2011 \\
\hline
\text{Pintz-Ruzsa}
&
8
&
2020 \\
\hline
\end{tabular}
\end{center}
\caption{Linnik-Goldbach representation}
\label{table8}
\end{table}

Before Linnik \cite{Linnik1953} proved the aforementioned result in 1953, he established a weaker conclusion in 1951. Under the assumption that the GRH is true, he proved that any sufficiently large even number can be expressed as the sum of two primes plus at most $K'$ powers of $2$ for some constant $K'$ \cite{Linnik1951}. This conclusion has also undergone a series of improvements. Heath-Brown and Puchta \cite{HBP2002} achieved the most refined estimate for $K'$ to date, surpassing the earlier results of Linnik \cite{Linnik1951}, Gallagher \cite{Gallagher1975}, Liu-Liu-Wang \cite{LLW1998_1}, Liu-Liu-Wang \cite{LLW1999}, and Wang \cite{Wang1999}.
\begin{theorem}[Heath-Brown--Puchta, 2002]
Suppose the GRH is true, every sufficiently large even integer is the sum of two primes and at most $7$ powers of $2$.
\end{theorem}

\begin{table}[ht]
\begin{center}
\begin{tabular}{ |C{6cm}|C{6cm}|C{1cm}| }
\hline
\text{Linnik}
&
$K'$
&
1951 \\
\hline
\text{Gallagher}
&
$K'$
&
1975 \\
\hline
\text{Liu-Liu-Wang}
&
770
&
1998 \\
\hline
\text{Liu-Liu-Wang}
&
200
&
1999 \\
\hline
\text{Wang}
&
160
&
1999 \\
\hline
\text{Heath-Brown--Puchta}
&
7
&
2002 \\
\hline
\end{tabular}
\end{center}
\caption{Linnik-Goldbach representation under GRH}
\label{table9}
\end{table}

\section{Some sufficient conditions for the Goldbach conjecture}\label{some sufficient conditions for the Goldbach conjecture}

In this section, we explore different approaches to establish sufficient conditions for the Goldbach conjecture to hold for sufficiently large even integers. This involves the use of the circle method, the study of the M\"{o}bius function, equi-distributions of certain arithmetic functions, and the behavior of the least prime in arithmetic progressions. These conditions not only provide a deeper understanding of the Goldbach conjecture but also offer potential avenues for further exploration.

\subsection{Goldbach conjecture and the circle method}

We have mentioned in Section~\ref{exceptionalsetofGoldbachnumberssec} that the circle method can be used to study the Goldbach conjecture and obtain information on the set of exceptional Goldbach numbers. To make the exceptional set of size $O(1)$, i.e., to prove the Goldbach conjecture for sufficiently large even integers, one needs to prove the integral over the so-called minor arc is smaller than that over the major arc.

We first introduce some notation. Let $N$ be a sufficiently large even integer and let $(\log N)^{15} \leq P \leq N^\epsilon$ for some $\epsilon > 0$. Let $x_0 = P/N$. A major arc $M(q, h)$ is defined to be the interval $[h/q - x_0, h/q + x_0]$, where $0 < h \leq q \leq P$ and $(h, q) = 1$. They are disjoint and contained in the closed interval $[x_0, 1 + x_0]$. The minor arc $m(N)$ is defined to be the set of points in $[x_0, 1 + x_0]$ which are not in any $M(q, h)$ for $(h, q) = 1$ and $q \leq P$, while $m^*(N)$ is defined to be the set of points in $[x_0, 1 + x_0]$ which are not in any $M(q, h)$ for $(h, q) = 1$, $(q, N) = 1$, and $q \leq P$. Then Mozzochi and Balasubramanian \cite{MB1978} proved in 1978 the following two theorems.
\begin{theorem}[Mozzochi-Balasubramanian, 1978]
Suppose the GRH is true, let $P = N^\epsilon$, if for any sufficiently large even integer $N$ we have
\[
\int_{m(N)} \left( \sum_{p \leq N} \exp(2 \pi i p x) \right)^2 \exp(- 2 \pi i N x) dx = o(N (\log N)^{-2}),
\]
then every sufficiently large even integer can be written as the sum of two primes.
\end{theorem}

\begin{theorem}[Mozzochi-Balasubramanian, 1978]
Let $P = N^{-1} (\log N)^{15}$, if for any sufficiently large even integer $N$ we have
\[
\int_{m^*(N)} \left( \sum_{p \leq N} \exp(2 \pi i p x) \right)^2 \exp(- 2 \pi i N x) dx = o(N (\log N)^{-2}),
\]
then every sufficiently large even integer can be written as the sum of two primes.
\end{theorem}

In 1983, they \cite{BM1983} proved two more related theorems, one of which involved weakening the assumption of the GRH to the assumption that Siegel zeros do not exist.

\begin{theorem}[Balasubramanian-Mozzochi, 1983]
Suppose Siegel zeros do not exist and let $P = \exp(c (\log N)^{1/2})$. If for any sufficiently large even integer $N$ we have
\[
\int_{m(N)} \left( \sum_{p \leq N} \exp(2 \pi i p x) \right)^2 \exp(- 2 \pi i N x) dx = o(N (\log N)^{-2}),
\]
then every sufficiently large even integer can be written as the sum of two primes.
\end{theorem}

\begin{theorem}[Balasubramanian-Mozzochi, 1983]
Let $P = \exp(c (\log N)^{1/2})$. If for any sufficiently large even integer $N$ we have
\[
\int_{m(N)} \left( \sum_{p \leq N} \exp(2 \pi i p x) \right)^2 \exp(- 2 \pi i N x) dx = o(N P^{- 1/32} (\log N)^{-2}),
\]
then every sufficiently large even integer can be written as the sum of two primes.
\end{theorem}

Mozzochi, in his research conducted in 1980 and 1981, delved deeper into these sufficient conditions for the Goldbach problem. He expressed these conditions in terms of integrals involving certain trigonometric functions \cite{Mozzochi1980, Mozzochi1981}. Furthermore, Mozzochi observed in 2014 that the sufficient condition for the twin prime conjecture, when expressed in a similar manner, is significantly simpler in comparison \cite{Mozzochi2014}.

\subsection{Goldbach conjecture and the M\"{o}bius function}

The M\"{o}bius function is defined as
\[
\mu(n) =
\begin{cases}
(-1)^{\nu(n)}, & \text{if } n \text{ is square-free}, \\
0, & \text{otherwise,}
\end{cases}
\]
where $\nu(n)$ is the number of distinct prime divisors of $n$. The M\"{o}bius function  plays an important role in number theory. For example, it can be used to express the principle of inclusion-exclusion, and there exists an equivalent form of the prime number theorem in terms of the M\"{o}bius function.

The Goldbach conjecture is also connected to the M\"{o}bius function. By studying the number of solutions of the Diophantine equation $n = ax + by$, Hua started an attempt work on the Goldbach conjecture since 1975 with the aid of Na and proved the following theorem \cite{Hua1989}.

\begin{theorem}[Hua, 1989]
For a sufficiently large even integer $N$, we let $\displaystyle H = \prod_{p \leq \sqrt{N}} p$ and let
\begin{align}\label{Huaeq}
\Phi(N) =
\sum_{d \mid (H, N)}
\sum_{\substack{a \left\vert \frac{H}{d} \right. \\ a \leq \frac{N}{d}}}
\sum_{\substack{b \left\vert \frac{H}{d} \right. \\ a \leq \frac{N}{d} \\ (a, b) = 1}} \mu(a) \mu(b)
\left(
\frac{1}{2} - \left\{ \frac{b^* N}{a} \right\}
\right),
\end{align}
where $b^* b \equiv 1 \pmod{a}$ and $\{x\}$ is the fractional part of $x$. If $\Phi(N) = o\left( \frac{N C_N}{(\log N)^2} \right)$, then every sufficiently large even integer can be written as the sum of two primes.
\end{theorem}

Na in 1986 \cite{Na1986} further studied the sum $\Phi(N)$ in \eqref{Huaeq}. He got rid of the congruence restriction $b^* b \equiv 1 \pmod{a}$ in~\eqref{Huaeq} and rewrote $\Phi(N)$ as
\begin{align*}
\Phi(N) & =
\sum_{d \mid (H, N)}
\sum_{\substack{a \left\vert \frac{H}{d} \right. \\ a \leq \frac{N}{d}}}
\sum_{\substack{b \left\vert \frac{H}{d} \right. \\ a \leq \frac{N}{d} \\ (a, b) = 1}} \mu(a) \mu(b)
\left(
\frac{1}{2} - \left\{ \frac{b^* N}{a} \right\}
\right) \\
& =
\sum_{d \mid (H, N)}
\sum_{\substack{a \left\vert \frac{H}{d} \right. \\ a \leq \frac{N}{d}}}
\sum_{\substack{b \left\vert \frac{H}{d} \right. \\ a \leq \frac{N}{d} \\ (a, b) = 1}}
\frac{\mu(a)}{a} \mu(b)
\sum_{r = 0}^{a - 1}
\left(
\frac{1}{2} - \frac{r}{a}
\right)
\sum_{\ell = 1}^a
e^{2 \pi i \left(\frac{N}{d} - b r\right) \ell/a}
\\
& =
\sum_{d \mid (H, N)}
\sum_{\substack{a \left\vert \frac{H}{d} \right. \\ a \leq \frac{N}{d}}}
\sum_{\substack{b \left\vert \frac{H}{d} \right. \\ a \leq \frac{N}{d} \\ (a, b) = 1}}
\mu(a) \mu(b)
\left(
\sum_{r = 0}^{a - 1}
\left(
\frac{1}{2} - \frac{r}{a}
\right)
+
\frac{1}{a}
\sum_{r = 0}^{a - 1}
\left(
\frac{1}{2} - \frac{r}{a}
\right)
\sum_{\ell = 1}^{a - 1}
e^{2 \pi i \left(\frac{N}{d} - b r\right) \ell/a}
\right)
\\
& =
\frac{1}{2}
\sum_{d \mid (H, N)}
\sum_{\substack{a \left\vert \frac{H}{d} \right. \\ a \leq \frac{N}{d}}}
\sum_{\substack{b \left\vert \frac{H}{d} \right. \\ a \leq \frac{N}{d} \\ (a, b) = 1}}
\frac{\mu(a)}{a} \mu(b)
-
\sum_{d \mid (H, N)}
\sum_{\substack{a \left\vert \frac{H}{d} \right. \\ a \leq \frac{N}{d}}}
\sum_{\substack{b \left\vert \frac{H}{d} \right. \\ a \leq \frac{N}{d} \\ (a, b) = 1}}
\frac{\mu(a)}{a^2} \mu(b)
\sum_{r = 1}^{a - 1}
\sum_{\ell = 1}^{a - 1}
r
e^{2 \pi i \left(\frac{N}{d} - b r\right) \ell/a}.
\end{align*}
Then he proved the first sum above equals $- \frac{N C_N}{2 (\log N)^3} + O \left( \frac{N \log \log N}{(\log N)^{7/2}} \right)$ and consequently obtained the following result.
\begin{theorem}[Na, 1986]
For a sufficiently large even integer $N$, we let $\displaystyle H = \prod_{p \leq \sqrt{N}} p$ and let
\[
\Phi'(N) =
\sum_{d \mid (H, N)}
\sum_{\substack{a \left\vert \frac{H}{d} \right. \\ a \leq \frac{N}{d}}}
\sum_{\substack{b \left\vert \frac{H}{d} \right. \\ a \leq \frac{N}{d} \\ (a, b) = 1}}
\frac{\mu(a)}{a^2} \mu(b)
\sum_{r = 1}^{a - 1}
\sum_{\ell = 1}^{a - 1}
r
e^{2 \pi i \left(\frac{N}{d} - b r\right) \ell/a}.
\]
If $\Phi'(N) = o\left( \frac{N C_N}{(\log N)^2} \right)$, then every sufficiently large even integer can be written as the sum of two primes.
\end{theorem}
We refer the readers to Baier's thesis \cite[Chapter 3]{Baier2000}
in 2000 for further study of Hua's approach to the Goldbach conjecture.

\subsection{Goldbach conjecture and equi-distributions of some arithmetic functions}

The prime number theorem $\pi(x) \sim \frac{x}{\log x}$ is usually stated in the form $\sum_{n \leq x} \Lambda(x) \sim x$, where $\Lambda(n)$ is the von Mangoldt function defined as
\[
\Lambda(n) =
\begin{cases}
\log p, & \text{if } n = p^\ell \text{ for some prime $p$ and integer $\ell \geq 1$}, \\
0, & \text{otherwise. }
\end{cases}
\]

Let $r(N)$ be the number of representations of $N$ as the sum of two primes, and  let $\widetilde{r}(N)$ denote the weighted number of representations of a large even integer $N$ as the sum of two primes
\begin{equation}\label{rneq}
\widetilde{r}(N)=\sum_{n < N} \Lambda(n) \Lambda(N - n).
\end{equation}
Since
\begin{equation*}
r(N) = \frac{\widetilde{r}(N)}{(\log N)^2}
\left(1 + O \left( \frac{\log \log N}{\log N} \right) \right) + O \left( \frac{N}{(\log N)^3} \right),
\end{equation*}
if we can obtain a good enough lower bound for $\widetilde{r}(N)$, it would imply that the Goldbach conjecture holds for sufficiently large even integers. Pan \cite{Pan1982} studied $\widetilde{r}(N)$ in 1982 and obtained the following result.

\begin{theorem}[Pan, 1982]
For sufficiently large even integers $N$ and for $Q = N^{\frac{1}{2}} (\log N)^{-20}$, we have the weighted number of representation of $N$ as the sum of two primes $\sum_{n < N} \Lambda(n) \Lambda(N - n)$ equals
\begin{align}\label{Paneq}
2C_N N
+ \sum_{n < N}
\left(
\sum_{\substack{d_1 \mid n \\ d_1 > Q}} \mu(d_1) \log d_1
\right)
\left(
\sum_{\substack{d_2 \mid N - n \\ (d_2, N) = 1 \\ d_2 > Q}} \mu(d_2) \log d_2
\right) + O\left( \frac{N}{\log N} \right).
\end{align}
If the triple sum in~\eqref{Paneq} is smaller than the expected main term $2C_N N$, then every sufficiently large even integer can be written as the sum of two primes.
\end{theorem}

Huang and Li \cite{HL2021} studied the triple sum in~\eqref{Paneq} and obtained the following theorem in 2021, which gives a sufficient condition for the Goldbach conjecture for sufficiently large even integers.
\begin{theorem}[Huang-L., 2021]
If the following inequality holds:
\[
\sum_{q \leq N^{\frac{1}{2} + \epsilon}}
\max_{y < N}
\max_{\substack{a \\ (a, q) = 1}}
\left|
\sum_{\substack{n \leq y \\ n \equiv a \bmod{q}}} \Lambda(N - n) \mu(n)
- \frac{1}{\phi(q)} \sum_{n \leq y} \Lambda(N - n) \mu(n)
\right|
\ll_A \frac{N}{(\log N)^A},
\]
then the Goldbach conjecture holds for all sufficiently large even integers.
\end{theorem}
We remark that an analogous result for the twin prime conjecture has been obtained by Ram Murty and Vatwani \cite{RV2017} in 2017. From these results we see it is important to study the behavior of the arithmetic functions $\Lambda(N - n) \mu(n)$ and $\Lambda(n - 2) \mu(n)$. Litchtman \cite{Lichtman2022} studied latter function on average in short intervals in 2022.

\subsection{Goldbach conjecture and the least prime in arithmetic progressions}

In 2010, Zhang made a innovative contribution to the study of the Goldbach conjecture  \cite{Zhang2010}, in which he found some connections between the Goldbach conjecture and the least prime in arithmetic progressions.

To state Zhang's result we first introduce some notation. For coprime positive integers $k$ and $\ell$ we let $p(k, \ell)$ be the least prime congruent to $\ell$ modulo $k$. We define $p(k)$ as the maximum value of $p(k, \ell)$, where $(k, \ell) = 1$ and $1 \leq \ell \leq k - 1$. In 2011, Xylouris proved that $p(k) \ll k^{5.18}$ \cite{Xylouris2011}. In 1934, Chowla had shown that $p(k) \ll k^{2 + \epsilon}$ for any $\epsilon > 0$ assuming the GRH \cite{Chowla1934} and he conjectured $p(k) \ll k^{1 + \epsilon}$.

We state Zhang's main theorem \cite{Zhang2010} as follows to end this section.
\begin{theorem}[Zhang, 2010]
If for every integer $n > 5$ there exists a natural number $r$ such that $2n - p_r$ is coprime to each of $2n - p_1, \dots, 2n - p_{r-1}, 2n - p_{r+1}, \dots, 2n - p_m$, where $p_1, \dots, p_{r-1}, p_r, p_{r+1}, \dots, p_m$ are all odd primes smaller than $n$, $p_r$ satisfies $(p_r, n) = 1$ and $1 \leq n \leq m = \pi(n - 1) - 1$, and for every sufficiently large positive integer $k$ we have $p(k) < k^{1.5}$, then every sufficiently large even integer is the sum of two primes.
\end{theorem}

\section{Some results of the same flavor as the Goldbach conjecture} \label{some results of the same flavor as the Goldbach conjecture}

In this section, we discuss some problems that share a similar nature with the Goldbach conjecture. The first type of problems involves representing natural numbers as the sum of integers with certain properties. Examples include representing positive odd integers as the sum of a prime and a power of $2$, representing positive integers as the sum of a prime and a square, representing positive integers as the sum of a prime and a Fibonacci number, and the Waring-Goldbach problem. Another type of questions explores the Goldbach problem in different algebraic structures, extending beyond the realm of integers. These variants include the Goldbach conjecture over number fields, the Goldbach conjecture over the polynomial ring with integral coefficients, the Goldbach conjecture over the formal power series ring with integral coefficients, the Goldbach conjecture over the polynomial ring with coefficients from a field, and the Goldbach conjecture for matrices, among others.

\subsection{Representing positive odd integers as the sum of a prime and a power of 2} De Polignac conjectured in 1849 that every odd number larger than $3$ can be written as the sum of an odd prime and a power of $2$. Soon he realized that this conjecture is false, as $127$ serves as a counterexample. In 1950, van der Corput \cite{vanderCorput1950} proved that the counterexamples have a positive density. Erd\H{o}s \cite{Erdos1950} proved the set of counterexamples contains an infinite arithmetic progression in the same year. Erd\H{o}s' result implies the upper density of numbers representable as the sum of a prime and a power of $2$ is at most $0.5 - 9 \times 10^{-8}$. This constant was improved to $0.4909$ by Habsieger and Roblot \cite{HR2006} in 2006.

\begin{table}[ht]
\begin{center}
\begin{tabular}{ |C{6cm}|C{6cm}|C{1cm}| }
\hline
\text{Erd\H{o}s}
&
0.49999991
&
1950 \\
\hline
\text{Habsieger-Roblot}
&
0.4909
&
2006 \\
\hline
\end{tabular}
\end{center}
\caption{Upper bound of Romanoff's constant}
\label{table10}
\end{table}

On the other hand, Romanoff \cite{Romanoff1934} proved in 1934 that such numbers have a positive lower density. Romanoff's constant has been improved by Chen-Sun \cite{CS2004}, L\"{u} \cite{Lv2007}, Habsieger-Roblot \cite{HR2006}, Pintz \cite{Pintz2006}, Habsieger-Sivak--Fischler \cite{HS2010}, and Elsholtz-Schlage--Puchta \cite{ES2018}.

\begin{table}[ht]
\begin{center}
\begin{tabular}{ |C{6cm}|C{6cm}|C{1cm}| }
\hline
\text{Romanoff}
&
$C > 0$
&
1934 \\
\hline
\text{Chen-Sun}
&
0.0868
&
2004 \\
\hline
\text{L\"{u}}
&
0.09322
&
2007 \\
\hline
\text{Habsieger-Roblot}
&
0.0933
&
2006 \\
\hline
\text{Pintz}
&
0.093626
&
2006 \\
\hline
\text{Habsieger-Sivak--Fischler}
&
0.0936275
&
2010 \\
\hline
\text{Elsholtz-Schlage--Puchta-Jan--Christoph}
&
0.107648
&
2018 \\
\hline
\end{tabular}
\end{center}
\caption{Lower bound of Romanoff's constant}
\label{table11}
\end{table}
By employing brute force computation and a probabilistic model by Bombieri, Romani \cite{Romani1983} made a conjecture in 1983 regarding the value of the Romanoff's constant, suggesting that the constant might be approximately equal to 0.434.

In his proof, Erd\H{o}s \cite{Erdos1950} used the covering system for $\mathbb{Z}$, which is a finite collection of arithmetic progressions whose union contains all integers. The covering system has become an interesting and useful tool in its own right.
Erd\H{o}s posed a number of problems concerning covering systems \cite{Erdos1957, Erdos1963, Erdos1971, Erdos1973, EG1980}, including the well-known minimum modulus problem: is there a uniform upper bound on the smallest modulus of all covering systems with distinct moduli? Building upon the work of Filaseta, Ford, Konyagin, Pomerance, and Yu \cite{FFSPY2007} as well as Nielsen \cite{Nielsen2009}, Hough \cite{Hough2015} resolved this problem in 2015, showing that the smallest modulus of all covering systems with distinct moduli is at most $10^{16}$. In 2022 Balister, Bollob\'as, Morris, Sahasrabudhe, and Tiba \cite{BBMRS2022} improved this bound to $616,000$. Additionally, Cummings, Filaseta, and Trifonov \cite{CFT2022} demonstrated in 2022 that the smallest modulus of all covering systems with distinct square-free moduli is at most $118$. Furthermore, Klein, Koukoulopoulos, and Lemieux \cite{KKL2022} in 2022 established an upper bound for the minimum modulus in covering systems with multiplicity $s$, where $s$ denotes the largest number of times a modulus appears in the system. Recently, Li, Wang, and Yi \cite{LWY2023} and Li, Wang, Wang, and Yi \cite{LWWY2023} generalized the result of Klein, Koukoulopoulos, and Lemieux in 2023 and proved a corresponding result in the number field setting and in the polynomial ring over the integers setting.

\subsection{Representing positive integers as the sum of a prime and a Fibonacci number}

Similar to de Polignac's question regarding representing positive integers as the sum of a prime and a power of $2$, an interesting topic is the representation of numbers as the sum of a prime and a Fibonacci number. Considering that the count of Fibonacci numbers and powers of $2$ up to a large positive number $X$ both have logarithmic growth $O(\log X)$, it is reasonable to guess that the results involving Romanoff's constant may extend to this type of representation as well. In 2010, Lee \cite{Lee2010} proved the following result.
\begin{theorem}[Lee, 2010]
The set of numbers that can be expressed as the sum of a prime and a Fibonacci number has a positive lower asymptotic density.
\end{theorem}

Liu and Xue \cite{LX2021} made Lee's result explicit in 2021.
\begin{theorem}[Liu-Xue, 2021]
The lower density of integers representable as the sum of a prime and a Fibonacci number is at least $0.0254905$.
\end{theorem}

Wang and Chen \cite{WC2023} improved the result by Liu and Xue by proving the following interesting result in 2023, in which they consider the integers that can be represented as the sum of a prime and a Fibonacci number in $1$ to $37$ ways.

\begin{theorem}[Wang-Chen, 2023]
The set of integers which can be represented as the sum of a prime and a Fibonacci number in at least $1$ and at most $37$ ways has the lower asymptotic density at least $0.143$.
\end{theorem}

\subsection{Representing positive integers as the sum of a prime and a square}

In 1923, Hardy and Littlewood proposed the following conjecture, among others, in their work \cite{HL1923}.
\begin{conjecture}[Hardy-Littlewood, 1923]
If $n$ is sufficiently large and not a square, then $n$ can be represented as the sum of a prime and a square.
\end{conjecture}

In 2003, Li \cite{Li2003primeandsquare} provided an estimate for the exceptional set of positive numbers in this representation, improving earlier works by Wang \cite{Wang1995}, Polyakov \cite{Polyakov1990}, Br\"{u}nner-Perelli-Pintz \cite{BPP1989}, Vinogradov \cite{Vinogradov1985}, and Davenport-Heilbronn \cite{DH1937}.
\begin{theorem}[Li, 2003]
Let $E(X)$ denote the number of natural numbers not exceeding $X$ which cannot be written as the sum of a prime and a square, then $E(X) \ll X^{0.982}$.
\end{theorem}

\begin{table}[ht]
\begin{center}
\begin{tabular}{ |C{6cm}|C{6cm}|C{1cm}| }
\hline
\text{Davenport-Heilbronn}
&
$E(X) \leq X (\log X)^{-C}$, \,\, $C > 0$
&
1937 \\
\hline
\text{Vinogradov}
&
$E(X) \leq X^\theta$, \,\, $\theta < 1$
&
1985 \\
\hline
\text{Br\"{u}nner-Perelli-Pintz}
&
$E(X) \leq X^\theta$, \,\, $\theta < 1$
&
1989 \\
\hline
\text{Polyakov}
&
$E(X) \leq X^\theta$, \,\, $\theta < 1$
&
1990 \\
\hline
\text{Wang}
&
$E(X) \leq X^{0.99}$
&
1995 \\
\hline
\text{Li}
&
$E(X) \leq X^{0.982}$
&
2003 \\
\hline
\end{tabular}
\end{center}
\caption{Sum of a prime and a square}
\label{table12}
\end{table}
Similarly, this representation has been studied under the assumption that the GRH is true and in short intervals \cite{Mikawa1993_2, PZ1995, PP1995, BP1996, Languasco2004, Brudern2008, LZ2008, Languasco2009, Suzuki2017}. The upper bound on the number of such representations is studied by Nayebi \cite{Nayebi2011} in 2011.

The representation of positive integers as the sum of a prime and twice a square is a lesser-known conjecture of Goldbach \cite{Hodges1993}: every odd integer could be written in the form $p + 2a^2$, where $p$ is a prime or is taken to be $1$, and $a \geq 0$ is an integer. Stein \cite{Stein1856} in 1856 found two counterexamples $5777$ and $5993$ to this conjecture by Goldbach. Guy \cite{Guy1990} pointed out in 1990 that the density of such counterexamples is zero. In fact, Hodges \cite{Hodges1993} made a stronger conjecture in 1993 that for every number $N \geq 1$ there are only a finite number of odd integers that cannot be represented as the sum of a prime and twice a square in at least $N$ ways.

Similarly, we can consider representing positive integers as the sum of a prime and a triangular number, or as the sum of a prime and a $k$-th power, or more generally, as the sum of a prime and an element from a set that contains values of a function evaluated at integers or just at primes. For example, Hardy and Littlewood conjectured in 1923 that \cite{HL1923} for integers $a > 0$ and $b$, let $N(n)$ be the number of representations of $n$ in the form $p + a m^2 + bm$, where $p$ is a prime and $m \in \mathbb{Z}$, then if $n$, $a$, and $b$ have a common factor, or if $n$ and $a + b$ are both even, or if $b^2 + 4 a n$ is a square, then $N(n) = o(\sqrt{n} / \log n)$, while in all other cases, we have
\[
N(n) \sim \frac{\epsilon}{\sqrt{a}} \frac{\sqrt{n}}{\log n}
\prod_{p \geq 2} \frac{p}{p - 1}
\prod_{\substack{p \nmid a \\ p \geq 3}} \left( 1 - \frac{1}{p - 1}
\legendre{b^2 + 4 a n}{p}
\right),
\]
where $\epsilon = 1$ if $a + b$ is odd and $\epsilon = 2$ if $a + b$ is even, and where $\legendre{b^2 + 4 a n}{p}$ is the Legendre's symbol.

Finally, we remark that there are other binary representations of large integers that involves primes in Beatty sequences \cite{Brudern2000, BG2017}, Piatetski-Shapiro primes \cite{Li2003hybrid, MW2006}, practical numbers \cite{PW2022}, Lehmer numbers \cite{SW2010}, and so on.

\subsection{Waring-Goldbach problem}

The Waring-Goldbach problem concerns the behaviour of the number of solutions $I_k(N)$ of the equation
\[
N = p_1^n + \cdots + p_k^n,
\]
where $p_1, \cdots, p_k$ are primes and $n \geq 1$. It has been proved \cite{Encyclopedia} that $I_k(N) > 0$ if $k = O(n \log n)$, while an asymptotic formula for $I_k(N)$ has been obtained for $k = O(n^2 \log n)$.

We refer the readers to Sections 5 and 18 of Vaughan's paper \cite{Vaughan2016} and Sections 1.3 and 1.4 of Kumchev and Tolev's paper \cite{KT2005} to learn more about the Waring-Goldbach problem and the problem of representing large integers as the sum of mixed powers.

\subsection{Finite Goldbach conjecture in number fields}

In 1956, Cohen \cite{Cohen1956} considered the finite Goldbach problem in algebraic number fields. Let $F$ be a number field of degree $n$ over $\mathbb{Q}$. Let $A$ be a proper ideal of $F$, and let $A = \mathfrak{P}_1^{\lambda_1} \cdots \mathfrak{P}_h^{\lambda_h} \mathfrak{Q}_1^{\mu_1} \cdots \mathfrak{Q}_k^{\mu_k}$ for distinct prime ideals $\mathfrak{P}_i$ and $\mathfrak{Q}_j$, where $\lambda_i > 0$, $N(\mathfrak{P}_i) > 2$ and $\mu_j > 0$, $N(\mathfrak{Q}_j) = 2$ for $1 \leq i \leq h$ and $1 \leq j \leq k$. Denote by $R(A)$ the ring of residue classes modulo $A$. One can define primes, composites, and unit residue classes in $R(A)$ as shown in \cite[Section 2]{Cohen1956}.

Building upon his earlier work \cite{Cohen1954} on the Goldbach properties in $\mathbb{Z} / m \mathbb{Z}$, Cohen \cite{Cohen1956} proved two theorems in 1956 that further expanded our understanding of Goldbach type properties within this algebraic structure.

\begin{theorem}[Cohen, 1956]
There exists an $s \geq 1$ such that all elements of $R(A)$ are expressible as sums of $s$ primes in $R(A)$ if and only if $h + k > 1$ and $h > 0$. For such ideals $A$, the minimum value $M$ of $s$ is given by $M = 2$ if $k = 0$ and $h \geq 2$; by $M = 3$ if (i) $k = 1$ and $h \geq 2$, if (ii) $h = k = \mu_1 = 1$, or (iii) if $k = 2$ and $h \geq 1$; by $M = 4$ if $h = 1$, $k = 1$, and $\mu_1 > 1$; and by $M = k$ if $k \geq 3$ and $h \geq 1$.
\end{theorem}

\begin{theorem}[Cohen, 1956]
There exists an $H \geq 1$ such that all elements of $R(A)$ are expressible as sums of at most $H$ primes in $R(A)$ if and only if $A$ satisfies neither $\ell = 1$ nor $h = 0$ while $k$ is odd. For all other proper ideals $A$, the number $H$ may be chosen to be $H = H'$, where $H' = 2$ if $k = 0$; $H' = 3$ if $h \geq 1$ and $k = 1$, $2$, or $3$; $H' = 4$ in case (1) $h \geq 1$ and $k = 4$, $5$, or $6$, and in case (ii) $h = 0$, $k = 2$; $H' = j + 1$ if $h \geq 1$, $k = 2j + 1 \geq 7$; $H' = j$ if $h \geq 1$, $k = 2j \geq 8$; and $H' = 2j$ if $h = 0$, $k = 2j \geq 4$. The minimum value $\theta$ of $H$ is given by $\theta = H'$ with these exceptions: $\theta = 2$ in case (i) $h = 1$, $k = 1$, $\lambda_1 = 1$, and in case (ii) $h = 0$, $k = 2$, $\mu_1 = 1$, $\mu_2 = 1$; $\theta = 3$ in case (i) $k = 2$, $h = 0$, and either $\mu_1 = 1$, $\mu_2 \neq 1$ or $\mu_1 \neq 1$, $\mu_2 = 1$, and in case (ii) $k = 4$, $h \geq 1$, and either $h \neq 1$ or $\mu_i = 1$ for at least one $i$.
\end{theorem}
Cohen also discussed the finite Goldbach property for algebraic number fields. We refer the readers to \cite[Section 4]{Cohen1956} for more details.

\subsection{Goldbach conjecture over the ring of integers of number fields}
The ring of integers $\mathcal{O}_K$ of an algebraic number field $K$ serves as the counterpart of $\mathbb{Z}$ in $\mathbb{Q}$. Mitsui's result in 1960 \cite{Mitsui1960} focuses on extending the Goldbach conjecture to the setting of $\mathcal{O}_K$, studying the representation of integers in $\mathcal{O}_K$ as sums of totally positive prime elements.
\begin{theorem}[Mitsui, 1960]
Let $K$ be an algebraic number field. Almost all totally positive even integers of $K$ can be represented as the sum of two totally positive odd primes of $K$.
\end{theorem}

Mitsui's approach to the Goldbach conjecture over the ring of integers of number fields is not the unique approach. Based on the calculus reformulation of the Goldbach conjecture that $g(x) = \left( \sum_p \frac{x^p}{p !} \right)$ has positive derivatives $g^{(2k)}(x)$ for every $k > 1$,  Knill \cite{Knill2016} in 2016 stated the Goldbach conjecture for Gaussian primes in $\mathcal{O}_{\mathbb{Q}(\sqrt{-1})}$ and Eisenstein primes in $\mathcal{O}_{\mathbb{Q}(\sqrt{-3})}$.

We remark that Knill \cite{Knill2016} also considered the Goldbach conjecture for Hurwitz primes in the set of quaternion integers $Q = \{(a, b, c, d) \mid a > 0, b > 0, c > 0, d > 0\}$ and Octavian primes in the division algebra $\mathbb{O} = \{(z, w) \mid z, w \in Q\}$.

\subsection{Goldbach conjecture over the polynomial ring with integral coefficients}

Hayes \cite{Hayes1965} proved the following result in 1965.
\begin{theorem}[Hayes, 1965]
In $\mathbb{Z}[x]$, every polynomial of degree $n \geq 1$ can be expressed as the sum of two irreducible polynomials each of degree $n$.
\end{theorem}
In 2006, Saidak \cite{Saidak2006} provided a concise proof of Hayes's result and extended it further by considering the number of representations of monic polynomials $f(x) \in \mathbb{Z}[x]$ as the sum of two irreducible monic polynomials $g(x)$ and $h(x)$ in $\mathbb{Z}[x]$. Specifically, for a given monic polynomial $f(x) \in \mathbb{Z}[x]$ of degree $d$, Saidak denoted $R(f, y)$ as the count of such representations where the coefficients of $g(x)$ and $h(x)$ are bounded by $y$ in absolute value. He proved that as $y$ approaches infinity, the number of representations $R(f, y)$ satisfies the inequality $y^d \ll_f R(f, y) \ll_f y^d$. In 2010, Kozek \cite{Kozek2007, Kozek2010} obtained the asymptotic formula for the number of such representations.

\begin{theorem}[Kozek, 2010]
Let $f(x)$ be a monic polynomial in $\mathbb{Z}[x]$ of degree $d \geq 1$. The number of representations of $f(x)$ as the sum of two irreducible monics $g(x)$ and $h(x)$ in $\mathbb{Z}[x]$, with the coefficients of $g(x)$ and $h(x)$ bounded in absolute value by $y$, is asymptotic to $(2 y)^{d - 1}$.
\end{theorem}
Additional generalizations and improvements have been made by Pollack \cite{Pollack2011} in 2011, Dubickas \cite{Dubikas2013} in 2013, and Lemos-de Araujo \cite{Ld2017} in 2017.

\subsection{Goldbach conjecture over the formal power series ring with integral coefficients}

Prime elements in the formal power series ring $\mathbb{Z}[[x]]$ have been very well characterized by Birmajer, Gil, and Weiner\cite{BGW2012} in 2012. Based on their work, Paran \cite{Paran2020} proved the following result in 2020.
\begin{theorem}[Paran, 2020]
A power series $f = f_0 + f_1 x + f_2 x^2 + \cdots \in \mathbb{Z}[[x]]$ is a sum of two primes in $\mathbb{Z}[[x]]$ if and only if $f_0$ is of the form $\pm p^k \pm q^\ell$ or of the form $\pm p^k$, where $p, q$ are prime integers and $k, \ell$ are positive integers.
\end{theorem}
It follows that not all elements of $\mathbb{Z}[[x]]$ may be written as a sum of two primes, since it has been shown by Sun \cite{Sun2000} that integers in the arithmetic progression
\[
f_0 \equiv 47867742232066880047611079 \pmod{66483034025018711639862527490}
\]
can not be written in the form $\pm p^k \pm q^\ell$ for any primes $p, q$ and infinitely many such elements are not prime powers.

\subsection{Goldbach conjecture over the polynomial ring with coefficients from a field}

For the ternary Goldbach problem in the polynomial ring over finite fields, Effinger and Hayes \cite{EH1991} proved the following result in 1991.

\begin{theorem}[Effinger-Hayes, 1991]
Every monic polynomial $f(x)$ of degree $r \geq 2$ over every finite field $\mathbb{F}_q$, except the cases $f(x) = x^2 + \alpha$ with $q$ even or $q = 2$ and for polynomials divisible by $x$ or $x + 1$, can be expressed as the sum of three irreducible monic polynomials in $\mathbb{F}_q[x]$, one of degree $r$ and the others of lower degree.
\end{theorem}

In 2014, Bender \cite{Bender2014} contributed to the binary Goldbach problem in the polynomial ring over finite fields and proved the following result.
\begin{theorem}[Bender, 2014]
Let $f(x) \in \mathbb{F}_q[x]$ be a monic polynomial of degree $n > 1$. If $q$ is odd and
\[
q > 3 n^4 \cdot 16^{n^4} \cdot (n + 1)^{8n^4 + 1},
\]
then it is the sum of two irreducible monic polynomials.
\end{theorem}

In 1998, Wang investigated the Goldbach problem in polynomial rings over infinite fields and obtained the following result.

\begin{theorem}[Wang, 1998]
There exist infinitely many infinite fields $F$ such that every monic polynomial of degree $d$ can be expressed as the sum of three irreducible monic polynomials in $F[x]$, one of degree $d$ and the others of lower degree, where $d \geq 2$ if $\text{char}(F)$ is 0 or odd and $d \geq 3$ if $\text{char}(F) = 2$.
\end{theorem}

\subsection{Goldbach conjecture over the ring of matrices}

A matrix in $M_n(\mathbb{Z})$ is irreducible if and only if its determinant equals $\pm p$ for some prime p. Thus, the corresponding Goldbach problem over $M_n(\mathbb{Z})$ can be stated as follows: can we write a matrix with even entries as the sum of two matrices, both of which have positive or negative prime determinants? Vaserstein \cite{Vaserstein1989} studied this problem for $M_2(\mathbb{Z})$ in 1989 and proved a stronger statement.
\begin{theorem}[Vaserstein, 1989]
Given any integer $p$ and any matrix $A \in M_2(\mathbb{Z})$, there are $X, Y \in M_2(\mathbb{Z})$ such that $A = X + Y$ and $\det(X) = \det(Y) = p$.
\end{theorem}
Vaserstein then raised a similar question for $M_3(\mathbb{Z})$. Wang \cite{Wang1992} provided an answer to this question in 1992 and derived more general results for $M_n(\mathbb{Z})$, where $n \geq 3$.

\begin{theorem}[Wang, 1992]
Vaserstein's result holds for $M_n(\mathbb{Z})$ for positive even integers $n$. Let $n > 1$ be an odd integer and $p$ a fixed integer. Then for any $A = (a_{ij})_{1 \leq i, j \leq n}$ in $M_n(\mathbb{Z})$ there are $X, Y \in M_n(\mathbb{Z})$ such that $A = X + Y$ and $\det(x) = \det(Y) = p$ if and only if $d(A) = 1$ or $2$, where $d(A) = \gcd\{a_{ij}\}_{1 \leq i, j \leq n}$.
\end{theorem}

Qin \cite{Qin1996} proved two related results in 1996.

\begin{theorem}[Qin, 1996]
For any $A \in M_n(\mathbb{Z})$ and any integer $p$, there are $X, Y \in M_n(\mathbb{Z})$ such that $A = X + Y$ and $\det(X) = (-1)^n \det(Y) = p$.
\end{theorem}

\begin{theorem}[Qin, 1996]
For any $A \in M_n(\mathbb{Z})$ and any integer $p$, if $n \geq 3$ is odd, then there are $X, Y, Z \in M_n(\mathbb{Z})$ such that $A = X + Y + Z$ and $\det(X) = \det(Y) = \det(Z) = p$.
\end{theorem}

In 1996, Bloy \cite{Bloy1996} proved the following result, which allows for a more general choice of determinants for the matrices $X$ and $Y$ in the decomposition.

\begin{theorem}[Bloy, 1996]
Fix $p$ and $q \in \mathbb{Z}$. There exist $X, Y \in M_n(\mathbb{Z})$ with $A = X + Y$, $\det(X) = p$, and $\det(Y) = q$ if and only if $d(A) \mid (p + (-1)^{n + 1} q)$.
\end{theorem}

In 2005, Hu \cite{Hu2005} generalized Bloy's results and proved the corresponding result for the matrix ring $M_n(R)$ over any principal ideal domain $R$.

Finally, we remark that in addition to his study of the Goldbach problem over $M_2(\mathbb{Z})$, Vaserstein also investigated several related problems in this setting. He considered Fermat's problem, Waring's problem, the commutators problem (i.e., the equation $XY = YX = A$ for a given matrix $A$), the problem of squares of linear forms and symmetric matrices, as well as primes of the form $X^2 + 1$ over $M_2(\mathbb{Z})$. Vaserstein's work has had a significant impact and has been followed up by many subsequent studies.

\section{Connections of the Goldbach conjecture to other conjectures} \label{connectionsec}

In this section, we explore the relationship between the Goldbach conjecture and two other famous conjectures: the twin prime conjecture and the Riemann hypothesis.

\subsection{Connection between the Goldbach conjecture and the twin prime conjecture}
Naturally, a comparison can be drawn between the Goldbach conjecture and the twin prime conjecture or the de Polignac's conjecture. The twin prime conjecture asserts the existence of infinitely many prime pairs differing by $2$, while de Polignac's conjecture extends this idea by proposing the presence of infinitely many prime pairs differing by any given positive even number. In this light, a more suitable parallel can be established between the Goldbach conjecture and the statement that for every positive even number, there exists at least one pair of primes differing by that precise value.

By simply modifying the sets used in sieving, Chen's theorem has a corresponding version related to the twin prime conjecture or the de Polignac's conjecture. In the same paper that Chen proved his famous $1+2$ theorem, Chen \cite{Chen1973} also proved the following result in 1973.

\begin{theorem}[Chen, 1973]
Every even integer can be expressed as the difference of infinitely many pairs of integers, where one integer is a prime and the other integer is the product of at most two primes.
\end{theorem}

Recently, there have been significant advancements in the twin prime conjecture. For instance, Zhang \cite{Zhang2014} proved the result that primes have finite gaps. Subsequently, Zhang's result of $70,000,000$ was improved to $246$ by the Polymath group \cite{Polymath2014}, showcasing the rapid progress in this area of research.
\begin{theorem}[Zhang, 2014]
There are infinitely many pairs of primes with gap at most $70$ million.
\end{theorem}

\begin{theorem}[Polymath, 2014, Theorem 4]
There are infinitely many pairs of primes with gap at most $246$.
\end{theorem}
Their proofs involve the concept of the \emph{admissible set}. For any $k \in \mathbb{N}$, an admissible $k$-tuple is a tuple $\mathcal{H} = (h_1, \dots, h_k)$ of $k$ increasing integers $h_1 < \dots < h_k$ which avoids at least one residue class modulo $p$ for every prime $p$. For any natural number $m$, let $H_m$ denote the quantity $H_m = \liminf_{n \to \infty} (p_{n + m} - p_n)$. The following theorem establishes a connection between the Goldbach conjecture and the twin prime conjecture, assuming the validity of the generalized Elliott-Halberstam conjecture (cf.~\cite[Claim~2.6]{Polymath2014}).

\begin{theorem}[Polymath, 2014, Theorem 5]\label{Polymath2014thm}
Assume the generalized Elliott-Halberstam conjecture GEH$[\theta]$ for all $0 < \theta < 1$. Then, at least one of the following statements is true.
\begin{itemize}
\item (Twin prime conjecture) $H_1 = 2$.
\item (near-miss to even Goldbach conjecture) If $n$ is a sufficiently large multiple of $6$, then at least one of $n$ and $n - 2$ is expressible as the sum of two primes, similarly with $n - 2$ replaced by $n + 2$. In particular, every sufficiently large even number lies within $2$ of the sum of two primes.
\end{itemize}
\end{theorem}

A result of this nature was obtained by Pintz \cite{Pintz2012} in 2012. He proved that either $H_1$ was finite or that every interval $[x, x + x^\epsilon]$ contained the sum of two primes if $x$ was sufficiently large depending on $\epsilon > 0$. Clearly this result of Pintz has been surpassed by Zhang's result on bounded gaps between primes. In 2021, Benatar \cite{Benatar2021} proved the following result of the same nature as Theorem~\ref{Polymath2014thm} and the aforementioned result by Pintz.

\begin{theorem}[Benatar, 2021]
Assuming that the primes have level of distribution greater than $1/2$, at least one of the following two properties holds.

\begin{itemize}
\item Consecutive Goldbach numbers lie within a finite distance from one another.

\item The set of de Polignac numbers has full density in the set of positive even numbers.
\end{itemize}
\end{theorem}
Finally, we remark that some of the methods employed to study the twin prime conjecture can also be applied to tackle the Goldbach conjecture. Usually the treatment of the Goldbach problem is trickier.

\subsection{Connection between the Goldbach conjecture and the Riemann hypothesis}

The Riemann Hypothesis states that the non-trivial zeros of the Riemann zeta function $\zeta(s)$ lie on the line with real part equal to 1/2. Granville \cite{Granville2007} proved the following result in 2007, which gives an equivalent form of the Riemann hypothesis in terms of the Goldbach representations.

\begin{theorem}
The Riemann hypothesis is equivalent to the estimate
\[
\sum_{2N \leq x} (G(2N) - J(2N)) \ll x^{3/2 + o(1)},
\]
where $G(2N) = \sum_{\substack{p + q = 2N \\ p, q \text{ prime}}} \log p \log q$ and $J(2N) = 2 C_2 N \prod_{\substack{p \mid N \\ p > 2}} \left( \frac{p - 1}{p - 2} \right)$.
\end{theorem}

Assuming the Riemann hypothesis, one may obtain asymptotic formulas for the average number representations of an even integer as the sum of two primes. For example, Languasco and Zaccagnini \cite{LZ2012, LZ2015} proved the following result in 2015.
\begin{theorem}[Languasco-Zaccagnini, 2015]
Assuming the Riemann hypothesis, we have
\[
\sum_{n \leq N}
\widetilde{r}(n)
= \frac{N^2}{2}
- \sum_\rho
\frac{N^{\rho + 1}}{\rho(\rho + 1)}
+ O(N (\log N)^3),
\]
where $\widetilde{r}(n)$ is defined in~\eqref{rneq} and the sum over $\rho$ runs over the complex zeros of the Riemann zeta-function.
\end{theorem}
We refer the readers to the paper by Goldston and Yang \cite{GY2016} in 2016 for more results of this type.

Wrapping up this section, we draw attention to Friedlander, Goldston, Iwaniec, and Surajaya's significant result \cite{FGIS2022} that establishes a connection between the Goldbach conjecture and the Siegel zeros in 2022, which is related to earlier works of Fei \cite{Fei2016}, Goldston-Suriajaya \cite{GS2021}, Bhowmik-Halupczok-Matsumoto-Suzuki \cite{BHMS2019}, Bhowmik-Halupczok \cite{BH2020}, and Jia \cite{Jia2022}.

\begin{theorem}[Friedlander-Goldston-Iwaniec-Surajaya, 2022]
Suppose for all sufficiently large even $N$ we have
\[
\delta
C_N N
< \sum_{\substack{m_1 + m_2 = N \\ 2 \nmid m_1 m_2}} \Lambda(m_1) \Lambda(m_2)
< (4 - \delta) C_N N
\]
for some fixed $0 < \delta < 2$, then there are no zeros $\rho = \sigma + i t$ of any Dirichlet $L$-function corresponding to character $\chi \pmod{q}$ in the region
\[
\sigma \geq 1 - c / \log q (|t| + 2)
\]
with a positive constant $c$ which is allowed to depend on $\delta$.
\end{theorem}

\section*{Acknowledgements}
The author's research is partially supported by the Fundamental Research Funds for the Central Universities, Nankai University (Grant No. 6323114) and the National Natural Science Foundation of China (Grant No. 12201313). The author would like to thank Runbo Li for some helpful discussions.

\bibliographystyle{plain}
\bibliography{bib}

\begin{thebibliography}{100}

\bibitem{Baier2000}
Stephan Baier.
\newblock \"{U}ber {Z}ug\"{a}nge von {H}ua und {P}an {C}hengdong zum
  {P}rimzahlzwillingsproblem.
\newblock {\em Dissertation am Fachbereich}, Mathematik/Informatik der FU
  Berlin, 2000.

\bibitem{Baker1967}
A.~Baker.
\newblock On some diophantine inequalities involving primes.
\newblock {\em J. Reine Angew. Math.}, 228:166--181, 1967.

\bibitem{BH1998}
R.~C. Baker and G.~Harman.
\newblock The three primes theorem with almost equal summands.
\newblock {\em R. Soc. Lond. Philos. Trans. Ser. A Math. Phys. Eng. Sci.},
  356(1738):763--780, 1998.

\bibitem{BM1983}
R.~Balasubramanian and C.~J. Mozzochi.
\newblock Siegel zeros and the {G}oldbach problem.
\newblock {\em J. Number Theory}, 16(3):311--332, 1983.

\bibitem{BBMRS2022}
Paul Balister, B\'{e}la Bollob\'{a}s, Robert Morris, Julian Sahasrabudhe, and
  Marius Tiba.
\newblock On the {E}rd{\H{o}}s covering problem: the density of the uncovered
  set.
\newblock {\em Invent. Math.}, 228(1):377--414, 2022.

\bibitem{BF1992}
Antal Balog and John Friedlander.
\newblock A hybrid of theorems of {V}inogradov and {P}iatetski-{S}hapiro.
\newblock {\em Pacific J. Math.}, 156(1):45--62, 1992.

\bibitem{BGN2007}
William~D. Banks, Ahmet~M. G\"{u}lo\u{g}lu, and C.~Wesley Nevans.
\newblock Representations of integers as sums of primes from a {B}eatty
  sequence.
\newblock {\em Acta Arith.}, 130(3):255--275, 2007.

\bibitem{Bauer2012}
Claus Bauer.
\newblock The binary {G}oldbach conjecture with restrictions on the primes.
\newblock {\em Far East J. Math. Sci. (FJMS)}, 70(1):87--120, 2012.

\bibitem{Bauer2017}
Claus Bauer.
\newblock Goldbach's conjecture in arithmetic progressions: number and size of
  exceptional prime moduli.
\newblock {\em Arch. Math. (Basel)}, 108(2):159--172, 2017.

\bibitem{Bauer2017_2}
Claus Bauer.
\newblock Large sieve inequality with sparse sets of moduli applied to
  {G}oldbach conjecture.
\newblock {\em Front. Math. China}, 12(2):261--280, 2017.

\bibitem{BW2013}
Claus Bauer and Yonghui Wang.
\newblock The binary {G}oldbach conjecture with primes in arithmetic
  progressions with large modulus.
\newblock {\em Acta Arith.}, 159(3):227--243, 2013.

\bibitem{BG2017}
A.~V. Begunts and D.~V. Goryashin.
\newblock Current problems connected to {B}eatty sequences.
\newblock {\em Chebyshevski\u{\i} Sb.}, 18(4):97--105, 2017.

\bibitem{Benatar2021}
Jacques {Benatar}.
\newblock {Goldbach versus de Polignac numbers}.
\newblock {\em arXiv e-prints}, page arXiv:1505.03104, May 2015.

\bibitem{Bender2014}
Andreas~O. Bender.
\newblock Representing an element in {$\mathbb{F}_q[t]$} as the sum of two
  irreducibles.
\newblock {\em Mathematika}, 60(1):166--182, 2014.

\bibitem{BH2020}
Gautami {Bhowmik} and Karin {Halupczok}.
\newblock {Condtional Bounds on Siegel Zeros}.
\newblock {\em arXiv e-prints}, page arXiv:2010.01308, October 2020.

\bibitem{BHMS2019}
Gautami Bhowmik, Karin Halupczok, Kohji Matsumoto, and Yuta Suzuki.
\newblock Goldbach representations in arithmetic progressions and zeros of
  {D}irichlet {$L$}-functions.
\newblock {\em Mathematika}, 65(1):57--97, 2019.

\bibitem{BGW2012}
Daniel Birmajer, Juan~B. Gil, and Michael Weiner.
\newblock Factoring polynomials in the ring of formal power series over
  {$\mathbb{Z}$}.
\newblock {\em Int. J. Number Theory}, 8(7):1763--1776, 2012.

\bibitem{Bloy1996}
Greg Bloy.
\newblock Goldbach's {P}roblem in {M}atrix {R}ings.
\newblock {\em Math. Mag.}, 69(2):136--137, 1996.

\bibitem{BD1966}
E.~Bombieri and H.~Davenport.
\newblock Small differences between prime numbers.
\newblock {\em Proc. Roy. Soc. London Ser. A}, 293:1--18, 1966.

\bibitem{Bordignon2022}
Matteo Bordignon.
\newblock An explicit version of {C}hen's theorem.
\newblock {\em Bull. Aust. Math. Soc.}, 105(2):344--346, 2022.

\bibitem{BJS2022}
Matteo {Bordignon}, Daniel~R. {Johnston}, and Valeriia {Starichkova}.
\newblock {An explicit version of Chen's theorem}.
\newblock {\em arXiv e-prints}, page arXiv:2207.09452, July 2022.

\bibitem{Brudern2000}
J\"{o}rg Br\"{u}dern.
\newblock Some additive problems of {G}oldbach's type.
\newblock {\em Funct. Approx. Comment. Math.}, 28:45--73, 2000.

\bibitem{Brudern2008}
J\"{o}rg Br\"{u}dern.
\newblock Representations of natural numbers as the sum of a prime and a
  {$k$}-th power.
\newblock {\em Tsukuba J. Math.}, 32(2):349--360, 2008.

\bibitem{BP1996}
J\"{o}rg Br\"{u}dern and Alberto Perelli.
\newblock The addition of primes and power.
\newblock {\em Canad. J. Math.}, 48(3):512--526, 1996.

\bibitem{Brun1920}
V.~Brun.
\newblock Le crible d'{E}ratosthene et le th\'{e}or\`{e}me de {G}oldbach.
\newblock {\em Selsk. Skr.}, 3, 1920.

\bibitem{BPP1989}
R.~Br\"{u}nner, A.~Perelli, and J.~Pintz.
\newblock The exceptional set for the sum of a prime and a square.
\newblock {\em Acta Math. Hungar.}, 53(3-4):347--365, 1989.

\bibitem{Cai2002}
Yingchun Cai.
\newblock Chen's theorem with small primes.
\newblock {\em Acta Math. Sin. (Engl. Ser.)}, 18(3):597--604, 2002.

\bibitem{Cai2008}
Yingchun Cai.
\newblock On {C}hen's theorem. {II}.
\newblock {\em J. Number Theory}, 128(5):1336--1357, 2008.

\bibitem{Cai2013}
Yingchun Cai.
\newblock A remark on the {G}oldbach-{V}inogradov theorem.
\newblock {\em Funct. Approx. Comment. Math.}, 48(part 1):123--131, 2013.

\bibitem{Cai2017}
Yingchun Cai.
\newblock Almost prime triples and {C}hen's theorem.
\newblock {\em Acta Arith.}, 179(3):233--250, 2017.

\bibitem{CL2011}
Yingchun Cai and Yingjie Li.
\newblock Chen's theorem with small primes.
\newblock {\em Chinese Ann. Math. Ser. B}, 32(3):387--396, 2011.

\bibitem{CL1999_2}
Yingchun Cai and Minggao Lu.
\newblock Chen's theorem in arithmetical progressions.
\newblock {\em Sci. China Ser. A}, 42(6):561--569, 1999.

\bibitem{CL1999}
Yingchun Cai and Minggao Lu.
\newblock Chen's theorem in short intervals.
\newblock {\em Acta Arith.}, 91(4):311--323, 1999.

\bibitem{CL2002}
Yingchun Cai and Minggao Lu.
\newblock On {C}hen's theorem.
\newblock In {\em Analytic number theory ({B}eijing/{K}yoto, 1999)}, volume~6
  of {\em Dev. Math.}, pages 99--119. Kluwer Acad. Publ., Dordrecht, 2002.

\bibitem{Chen1965}
Jingrun Chen.
\newblock On large odd numbers as sum of three almost equal primes.
\newblock {\em Sci. Sinica}, 14:1113--1117, 1965.

\bibitem{Chen1966}
Jingrun Chen.
\newblock On the representation of a large even integer as the sum of a prime
  and the product of at most two primes.
\newblock {\em Kexue Tongbao}, 17:385--386, 1966.

\bibitem{Chen1973}
Jingrun Chen.
\newblock On the representation of a larger even integer as the sum of a prime
  and the product of at most two primes.
\newblock {\em Sci. Sinica}, 16:157--176, 1973.

\bibitem{Chen1978_2}
Jingrun Chen.
\newblock Further improvement on the constant in the proposition `1+2': On the
  representation of a large even integer as the sum of a prime and the product
  of at most two primes (ii).
\newblock {\em Sci. Sinica}, pages 477--494(in Chinese), 1978.

\bibitem{Chen1978_3}
Jingrun Chen.
\newblock On the {G}oldbach's problem and the sieve methods.
\newblock {\em Sci. Sinica}, 21(6):701--739, 1978.

\bibitem{Chen1978_1}
Jingrun Chen.
\newblock On the representation of a large even integer as the sum of a prime
  and the product of at most two primes. {II}.
\newblock {\em Sci. Sinica}, 21(4):421--430, 1978.

\bibitem{Chen1983}
Jingrun Chen.
\newblock The exceptional set of {G}oldbach numbers. {II}.
\newblock {\em Sci. Sinica Ser. A}, 26(7):714--731, 1983.

\bibitem{CL1989}
Jingrun Chen and Jianmin Liu.
\newblock The exceptional set of {G}oldbach-numbers. {III}.
\newblock {\em Chinese Quart. J. Math.}, 4(1):1--15, 1989.

\bibitem{CP1980}
Jingrun Chen and Chengdong Pan.
\newblock The exceptional set of {G}oldbach numbers. {I}.
\newblock {\em Sci. Sinica}, 23(4):416--430, 1980.

\bibitem{CS2004}
Yonggao Chen and Xuegong Sun.
\newblock On {R}omanoff's constant.
\newblock {\em J. Number Theory}, 106(2):275--284, 2004.

\bibitem{CT2017}
Tak~Wing Ching and Kaiman Tsang.
\newblock Small prime solutions to linear equations in three variables.
\newblock {\em Acta Arith.}, 178(1):57--76, 2017.

\bibitem{Choi1997}
Kwok Kwong~Stephen Choi.
\newblock A numerical bound for {B}aker's constant---some explicit estimates
  for small prime solutions of linear equations.
\newblock {\em Bull. Hong Kong Math. Soc.}, 1(1):1--19, 1997.

\bibitem{CK2006}
Kwok Kwong~Stephen Choi and Angel~V. Kumchev.
\newblock Mean values of {D}irichlet polynomials and applications to linear
  equations with prime variables.
\newblock {\em Acta Arith.}, 123(2):125--142, 2006.

\bibitem{Chowla1934}
S.~Chowla.
\newblock On the least prime in an arithmetic progression.
\newblock {\em J. Indian Math. Soc.}, 1(2):1--3, 1934.

\bibitem{Chudakov1937}
N.~G. Chudakov.
\newblock Sur le probl$\acute{e}$me de {G}oldbach.
\newblock {\em C.R.(Dokl.) Acad. Sci. URSS}, 17:335--338, 1937.

\bibitem{Cohen1954}
Eckford Cohen.
\newblock A finite analogue of the {G}oldbach problem.
\newblock {\em Proc. Amer. Math. Soc.}, 5:478--483, 1954.

\bibitem{Cohen1956}
Eckford Cohen.
\newblock The finite {G}oldbach problem in algebraic number fields.
\newblock {\em Proc. Amer. Math. Soc.}, 7:500--506, 1956.

\bibitem{vanderCorput1937}
J.~G. van~der Corput.
\newblock Sur l'hypoth$\acute{e}$se de {G}oldbach pour {P}resque tous les
  nombres pairs.
\newblock {\em Acta Arith}, 2:266--290, 1937.

\bibitem{vanderCorput1950}
J.~G. van~der Corput.
\newblock On de {P}olignac's conjecture.
\newblock {\em Simon Stevin}, 27:99--105, 1950.

\bibitem{CFT2022}
Maria {Cummings}, Michael {Filaseta}, and Ognian {Trifonov}.
\newblock {An upper bound for the minimum modulus in a covering system with
  squarefree moduli}.
\newblock {\em arXiv e-prints}, page arXiv:2211.08548, November 2022.

\bibitem{DH1937}
H.~Davenport and H.~Heilbronn.
\newblock Note on a {R}esult in the {A}dditive {T}heory of {N}umbers.
\newblock {\em Proc. London Math. Soc. (2)}, 43(2):142--151, 1937.

\bibitem{Davenport2000}
Harold Davenport.
\newblock {\em Multiplicative number theory}, volume~74 of {\em Graduate Texts
  in Mathematics}.
\newblock Springer-Verlag, New York, third edition, 2000.

\bibitem{Dubikas2013}
Art\={u}ras Dubickas.
\newblock Linear forms in monic integer polynomials.
\newblock {\em Canad. Math. Bull.}, 56(3):510--519, 2013.

\bibitem{Dudek2017}
A.~W. Dudek.
\newblock On the sum of a prime and a square-free number.
\newblock {\em Ramanujan J.}, 42(1):233--240, 2017.

\bibitem{Dufner1994}
Gunter Dufner.
\newblock Bin{\"a}res {G}oldbachproblem in kurzen {I}ntervallen teil i: Die
  explizite formel.
\newblock {\em Periodica Mathematica Hungarica}, 29:213--243, 1994.

\bibitem{Dufner1995}
Gunter Dufner.
\newblock Bin{\"a}res {G}oldbachproblem in kurzen {I}ntervallen ii.
\newblock {\em Periodica Mathematica Hungarica}, 30:37--60, 1995.

\bibitem{EH1991}
Gove~W. Effinger and David~R. Hayes.
\newblock A complete solution to the polynomial {$3$}-primes problem.
\newblock {\em Bull. Amer. Math. Soc. (N.S.)}, 24(2):363--369, 1991.

\bibitem{ES2018}
Christian Elsholtz and Jan-Christoph Schlage-Puchta.
\newblock On {R}omanov's constant.
\newblock {\em Math. Z.}, 288(3-4):713--724, 2018.

\bibitem{Erdos1950}
Paul Erd\H{o}s.
\newblock On integers of the form {$2^k+p$} and some related problems.
\newblock {\em Summa Brasil. Math.}, 2:113--123, 1950.

\bibitem{Erdos1957}
Paul Erd\H{o}s.
\newblock Some unsolved problems.
\newblock {\em Michigan Math. J.}, 4:291--300, 1957.

\bibitem{Erdos1963}
Paul Erd\H{o}s.
\newblock Quelques probl\`emes de th\'{e}orie des nombres.
\newblock In {\em Monographies de {L}'{E}nseignement {M}ath\'{e}matique, {N}o.
  6}, pages 81--135. Universit\'{e} de Gen\`eve, L'Enseignement
  Math\'{e}matique, Geneva, 1963.

\bibitem{Erdos1971}
Paul Erd\H{o}s.
\newblock Some problems in number theory.
\newblock Computers in {Number} {Theory}, {Proc}. {Atlas} {Sympos}. {No}. 2,
  {Oxford} 1969, 405-414 (1971)., 1971.

\bibitem{Erdos1973}
Paul Erd\H{o}s.
\newblock R\'{e}sultats et probl\`emes en th\'{e}orie des nombres.
\newblock In {\em S\'{e}minaire {D}elange-{P}isot-{P}oitou (14e ann\'{e}e:
  1972/73), {T}h\'{e}orie des nombres, {F}asc. 2}, pages Exp. No. 24, 7.
  Secr\'{e}tariat Math\'{e}matique, Paris, 1973.

\bibitem{EG1980}
Paul Erd\H{o}s and Ronald~Lewis Graham.
\newblock {\em Old and new problems and results in combinatorial number
  theory}, volume~28 of {\em Monographies de L'Enseignement Math\'{e}matique
  [Monographs of L'Enseignement Math\'{e}matique]}.
\newblock Universit\'{e} de Gen\`eve, L'Enseignement Math\'{e}matique, Geneva,
  1980.

\bibitem{Estermann1931}
T.~Estermann.
\newblock On the representations of a number as the sum of a prime and a
  quadratfrei number.
\newblock {\em Journal of the London Mathematical Society}, s1-6(3):219--221,
  1931.

\bibitem{Estermann1938}
T.~Estermann.
\newblock On {G}oldbach's {P}roblem : {P}roof that {A}lmost all {E}ven
  {P}ositive {I}ntegers are {S}ums of {T}wo {P}rimes.
\newblock {\em Proc. London Math. Soc. (2)}, 44(4):307--314, 1938.

\bibitem{Fei2016}
Jinhua Fei.
\newblock An application of the {H}ardy-{L}ittlewood conjecture.
\newblock {\em J. Number Theory}, 168:39--44, 2016.

\bibitem{FFSPY2007}
Michael Filaseta, Kevin Ford, Sergei Konyagin, Carl Pomerance, and Gang Yu.
\newblock Sieving by large integers and covering systems of congruences.
\newblock {\em J. Amer. Math. Soc.}, 20(2):495--517, 2007.

\bibitem{FrancisLee2020}
F.~J. {Francis} and E.~S. {Lee}.
\newblock {Additive Representations of Natural Numbers}.
\newblock {\em arXiv e-prints}, page arXiv:2003.08083, March 2020.

\bibitem{FrancisLee2022}
F.~J. {Francis} and E.~S. {Lee}.
\newblock {Additive Representations of Natural Numbers}.
\newblock {\em Integers}, 22:Paper No. A14, 10, 2022.

\bibitem{FKS2021}
C.~Frei, P.~Koymans, and E.~Sofos.
\newblock Vinogradov's three primes theorem with primes having given primitive
  roots.
\newblock {\em Math. Proc. Cambridge Philos. Soc.}, 170(1):75--110, 2021.

\bibitem{FGIS2022}
J.~B. Friedlander, D.~A. Goldston, H.~Iwaniec, and A.~I. Suriajaya.
\newblock Exceptional zeros and the {G}oldbach problem.
\newblock {\em J. Number Theory}, 233:78--86, 2022.

\bibitem{Gallagher1975}
P.~X. Gallagher.
\newblock Primes and powers of {$2$}.
\newblock {\em Invent. Math.}, 29(2):125--142, 1975.

\bibitem{GS2021}
D.~A. {Goldston} and Ade~Irma {Suriajaya}.
\newblock {Note on the Goldbach Conjecture and Landau-Siegel Zeros}.
\newblock {\em arXiv e-prints}, page arXiv:2104.09407, April 2021.

\bibitem{GY2016}
D.~A. Goldston and Liyang Yang.
\newblock {The Average Number of Goldbach Representations}.
\newblock {\em arXiv e-prints}, page arXiv:1601.06902, January 2016.

\bibitem{GVT1988}
A.~Granville, J.~van~de Lune, and H.~J.~J. te~Riele.
\newblock Checking the {G}oldbach conjecture on a vector computer.
\newblock In {\em Number theory and applications ({B}anff, {AB}, 1988)}, volume
  265 of {\em NATO Adv. Sci. Inst. Ser. C: Math. Phys. Sci.}, pages 423--433.
  Kluwer Acad. Publ., Dordrecht, 1989.

\bibitem{Granville2007}
Andrew Granville.
\newblock Refinements of {G}oldbach's conjecture, and the generalized {R}iemann
  hypothesis.
\newblock {\em Funct. Approx. Comment. Math.}, 37(part 1):159--173, 2007.

\bibitem{Grimmelt2022}
Lasse Grimmelt.
\newblock Goldbach numbers in short intervals.
\newblock {\em Ann. Sc. Norm. Super. Pisa Cl. Sci. (5)}, 23(3):1395--1416,
  2022.

\bibitem{Grimmelt2022FIprimes}
Lasse Grimmelt.
\newblock Vinogradov's theorem with {F}ouvry-{I}waniec primes.
\newblock {\em Algebra Number Theory}, 16(7):1705--1776, 2022.

\bibitem{GT2022}
Lasse {Grimmelt} and Joni {Ter{\"a}v{\"a}inen}.
\newblock {The Exceptional Set in Goldbach's Problem with Almost Twin Primes}.
\newblock {\em arXiv e-prints}, page arXiv:2207.08805, July 2022.

\bibitem{Guy1990}
Richard~K. Guy.
\newblock The second strong law of small numbers.
\newblock {\em Math. Mag.}, 63(1):3--20, 1990.

\bibitem{HR2006}
Laurent Habsieger and Xavier-Fran\c{c}ois Roblot.
\newblock On integers of the form {$p+2^k$}.
\newblock {\em Acta Arith.}, 122(1):45--50, 2006.

\bibitem{HS2010}
Laurent Habsieger and Jimena Sivak-Fischler.
\newblock An effective version of the {B}ombieri-{V}inogradov theorem, and
  applications to {C}hen's theorem and to sums of primes and powers of two.
\newblock {\em Arch. Math. (Basel)}, 95(6):557--566, 2010.

\bibitem{Halberstam1975}
Heini Halberstam.
\newblock A proof of {C}hen's theorem.
\newblock In {\em Journ\'{e}es {A}rithm\'{e}tiques de {B}ordeaux ({C}onf.,
  {U}niv. {B}ordeaux, {B}ordeaux, 1974)}, volume No. 24-25 of {\em
  Ast\'{e}risque}, pages 281--293. Soc. Math. France, Paris, 1975.

\bibitem{HL1923}
G.~H. Hardy and J.~E. Littlewood.
\newblock Some problems of `{P}artitio numerorum'; {III}: {O}n the expression
  of a number as a sum of primes.
\newblock {\em Acta Math.}, 44(1):1--70, 1923.

\bibitem{HL1924}
G.~H. Hardy and J.~E. Littlewood.
\newblock Some problems of '{P}artitio numerorum'({V}): {A} further
  contribution to the study of {G}oldbach's problem.
\newblock {\em Proc. London Math. Soc. (2)}, 22:46--56, 1924.

\bibitem{Harman2007}
Glyn Harman.
\newblock {\em Prime-detecting sieves}, volume~33 of {\em London Mathematical
  Society Monographs Series}.
\newblock Princeton University Press, Princeton, NJ, 2007.

\bibitem{Harman2020}
Glyn Harman.
\newblock Diophantine approximation with {G}oldbach numbers.
\newblock {\em Funct. Approx. Comment. Math.}, 63(2):151--163, 2020.

\bibitem{Haselgrove1951}
C.~B. Haselgrove.
\newblock Some theorems in the analytic theory of numbers.
\newblock {\em J. London Math. Soc.}, 26:273--277, 1951.

\bibitem{HJ2021}
S.~{Hathi} and D.~R. {Johnston}.
\newblock {On the sum of a prime and a square-free number with divisibility
  conditions}.
\newblock {\em arXiv e-prints}, page arXiv:2109.11883, September 2021.

\bibitem{Hayes1965}
D.~R. Hayes.
\newblock A {G}oldbach theorem for polynomials with integral coefficients.
\newblock {\em Amer. Math. Monthly}, 72:45--46, 1965.

\bibitem{HBP2002}
D.~R. Heath-Brown and J.-C. Puchta.
\newblock Integers represented as a sum of primes and powers of two.
\newblock {\em Asian J. Math.}, 6(3):535--565, 2002.

\bibitem{HBL2016}
Roger Heath-Brown and Xiannan Li.
\newblock Almost prime triples and {C}hen's theorem.
\newblock {\em J. Number Theory}, 169:265--294, 2016.

\bibitem{Helfgott2012Minor}
H.~A. {Helfgott}.
\newblock {Minor arcs for Goldbach's problem}.
\newblock {\em arXiv e-prints}, page arXiv:1205.5252, May 2012.

\bibitem{Helfgott2013Major}
H.~A. {Helfgott}.
\newblock {Major arcs for Goldbach's problem}.
\newblock {\em arXiv e-prints}, page arXiv:1305.2897, May 2013.

\bibitem{Helfgott2013}
H.~A. {Helfgott}.
\newblock {The ternary Goldbach conjecture is true}.
\newblock {\em arXiv e-prints}, page arXiv:1312.7748, December 2013.

\bibitem{Helfgott2015}
H.~A. {Helfgott}.
\newblock {The ternary Goldbach problem}.
\newblock {\em arXiv e-prints}, page arXiv:1501.05438, January 2015.

\bibitem{Hodges1993}
Laurent Hodges.
\newblock A {L}esser-{K}nown {G}oldbach {C}onjecture.
\newblock {\em Math. Mag.}, 66(1):45--47, 1993.

\bibitem{Hough2015}
Bob Hough.
\newblock Solution of the minimum modulus problem for covering systems.
\newblock {\em Ann. of Math. (2)}, 181(1):361--382, 2015.

\bibitem{Hu2005}
Wei Hu.
\newblock Goldbach's problem in the matrix ring over a principal ideal domain.
\newblock {\em Northeast. Math. J.}, 21(3):355--364, 2005.

\bibitem{Hua1989}
Lookeng Hua.
\newblock A direct attempt to {G}oldbach problem.
\newblock {\em Acta Math. Sinica (N.S.)}, 5(1):1--8, 1989.

\bibitem{HL2021}
Jingjing Huang and Huixi Li.
\newblock On the connection between the {G}oldbach conjecture and the
  {E}lliott-{H}alberstam conjecture.
\newblock In {\em Combinatorial and additive number theory {IV}}, volume 347 of
  {\em Springer Proc. Math. Stat.}, pages 323--346. Springer, Cham, 2021.

\bibitem{Jia1989}
Chaohua Jia.
\newblock The three-primes theorem over short intervals.
\newblock {\em Acta Math. Sinica}, 32(4):464--473, 1989.

\bibitem{Jia1991_1}
Chaohua Jia.
\newblock Three primes theorem in a short interval. {II}.
\newblock In {\em International {S}ymposium in {M}emory of {H}ua {L}oo {K}eng,
  {V}ol. {I} ({B}eijing, 1988)}, pages 103--115. Springer, Berlin, 1991.

\bibitem{Jia1991_2}
Chaohua Jia.
\newblock Three primes theorem in a short interval. {III}.
\newblock {\em Sci. China Ser. A}, 34(9):1039--1056, 1991.

\bibitem{Jia1991_3}
Chaohua Jia.
\newblock Three primes theorem in a short interval. {IV}.
\newblock {\em Adv. in Math. (China)}, 20(1):109--126, 1991.

\bibitem{Jia1992}
Chaohua Jia.
\newblock Three primes theorem in a short interval. {V}.
\newblock {\em Acta Math. Sinica (N.S.)}, 7(2):135--170, 1991.

\bibitem{Jia1991_4}
Chaohua Jia.
\newblock Three primes theorem in a short interval. {VI}.
\newblock {\em Acta Math. Sinica}, 34(6):832--850, 1991.

\bibitem{Jia1994}
Chaohua Jia.
\newblock Three primes theorem in a short interval. {VII}.
\newblock {\em Acta Math. Sinica (N.S.)}, 10(4):369--387, 1994.

\bibitem{Jia1995short}
Chaohua Jia.
\newblock Goldbach numbers in short interval. {I}.
\newblock {\em Sci. China Ser. A}, 38(4):385--406, 1995.

\bibitem{Jia1995short2}
Chaohua Jia.
\newblock Goldbach numbers in short interval. {II}.
\newblock {\em Sci. China Ser. A}, 38(5):513--523, 1995.

\bibitem{Jia1995}
Chaohua Jia.
\newblock On the {P}iatetski-{S}hapiro-{V}inogradov theorem.
\newblock {\em Acta Arith.}, 73(1):1--28, 1995.

\bibitem{Jia1996}
Chaohua Jia.
\newblock Almost all short intervals containing prime numbers.
\newblock {\em Acta Arith.}, 76(1):21--84, 1996.

\bibitem{Jia1996short}
Chaohua Jia.
\newblock On the exceptional set of {G}oldbach numbers in a short interval.
\newblock {\em Acta Arith.}, 77(3):207--287, 1996.

\bibitem{Jia2022}
Chaohua Jia.
\newblock On the conditional bounds for {S}iegel zeros.
\newblock {\em Acta Math. Sin. (Engl. Ser.)}, 38(5):869--876, 2022.

\bibitem{Kan1991}
Jiahai Kan.
\newblock On the number of solutions of {$N-p=P_r$}.
\newblock {\em J. Reine Angew. Math.}, 414:117--130, 1991.

\bibitem{Kan1992}
Jiahai Kan.
\newblock On the problem of {G}oldbach's type.
\newblock {\em Math. Ann.}, 292(1):31--42, 1992.

\bibitem{KKL2022}
Jonah {Klein}, Dimitris {Koukoulopoulos}, and Simon {Lemieux}.
\newblock {On the $j$-th smallest modulus of a covering system with distinct
  moduli}.
\newblock {\em arXiv e-prints}, page arXiv:2212.01299, December 2022.

\bibitem{Knill2016}
Oliver {Knill}.
\newblock {Goldbach for Gaussian, Hurwitz, Octavian and Eisenstein primes}.
\newblock {\em arXiv e-prints}, page arXiv:1606.05958, June 2016.

\bibitem{Kozek2010}
Mark Kozek.
\newblock An asymptotic formula for {G}oldbach's conjecture with monic
  polynomials in {$\mathbb{Z}[x]$}.
\newblock {\em Amer. Math. Monthly}, 117(4):365--369, 2010.

\bibitem{Kozek2007}
Mark~Robert Kozek.
\newblock {\em Applications of covering systems of integers and {G}oldbach's
  conjecture for monic polynomials}.
\newblock ProQuest LLC, Ann Arbor, MI, 2007.

\bibitem{Kumchev1997}
Angel~V. Kumchev.
\newblock On the {P}iatetski-{S}hapiro-{V}inogradov theorem.
\newblock {\em J. Th\'{e}or. Nombres Bordeaux}, 9(1):11--23, 1997.

\bibitem{Kumchev2008}
Angel~V. Kumchev.
\newblock On sums of primes from {B}eatty sequences.
\newblock {\em Integers}, 8:A8, 12, 2008.

\bibitem{KT2005}
Angel~V. Kumchev and D.~I. Tolev.
\newblock An invitation to additive prime number theory.
\newblock {\em Serdica Math. J.}, 31(1-2):1--74, 2005.

\bibitem{Languasco2004Goldbach}
Alessandro Languasco.
\newblock On the exceptional set of {G}oldbach's problem in short intervals.
\newblock {\em Monatsh. Math.}, 141(2):147--169, 2004.

\bibitem{Languasco2004}
Alessandro Languasco.
\newblock On the exceptional set of {H}ardy-{L}ittlewood's numbers in short
  intervals.
\newblock {\em Tsukuba J. Math.}, 28(1):169--191, 2004.

\bibitem{Languasco2009}
Alessandro Languasco.
\newblock A conditional result on the exceptional set for {H}ardy-{L}ittlewood
  numbers in short intervals.
\newblock {\em Int. J. Number Theory}, 5(6):933--951, 2009.

\bibitem{LZ2008}
Alessandro Languasco and Alessandro Zaccagnini.
\newblock On the {H}ardy-{L}ittlewood problem in short intervals.
\newblock {\em Int. J. Number Theory}, 4(5):715--723, 2008.

\bibitem{LZ2012}
Alessandro Languasco and Alessandro Zaccagnini.
\newblock The number of {G}oldbach representations of an integer.
\newblock {\em Proc. Amer. Math. Soc.}, 140(3):795--804, 2012.

\bibitem{LZ2015}
Alessandro Languasco and Alessandro Zaccagnini.
\newblock A {C}es\`aro average of {G}oldbach numbers.
\newblock {\em Forum Math.}, 27(4):1945--1960, 2015.

\bibitem{Lee2010}
K.~S.~Enoch Lee.
\newblock On the sum of a prime and a {F}ibonacci number.
\newblock {\em Int. J. Number Theory}, 6(7):1669--1676, 2010.

\bibitem{Ld2017}
Ab\'{\i}lio Lemos and Anderson Lu\'{\i}s~Albuquerque de~Araujo.
\newblock An asymptotic formula for {G}oldbach's conjecture with monic
  polynomials in {$\mathbb{Z}[\theta][x]$}.
\newblock {\em Colloq. Math.}, 148(2):215--223, 2017.

\bibitem{Li1995}
Hongze Li.
\newblock Goldbach numbers in short intervals.
\newblock {\em Sci. China Ser. A}, 38(6):641--652, 1995.

\bibitem{Li1999}
Hongze Li.
\newblock The exceptional set of {G}oldbach numbers.
\newblock {\em Quart. J. Math. Oxford Ser. (2)}, 50(200):471--482, 1999.

\bibitem{Li2000exp}
Hongze Li.
\newblock The exceptional set of {G}oldbach numbers. {II}.
\newblock {\em Acta Arith.}, 92(1):71--88, 2000.

\bibitem{Li2000Linnik}
Hongze Li.
\newblock The number of powers of 2 in a representation of large even integers
  by sums of such powers and of two primes.
\newblock {\em Acta Arith.}, 92(3):229--237, 2000.

\bibitem{Li2001}
Hongze Li.
\newblock The number of powers of 2 in a representation of large even integers
  by sums of such powers and of two primes. {II}.
\newblock {\em Acta Arith.}, 96(4):369--379, 2001.

\bibitem{Li2001smallprime}
Hongze Li.
\newblock Small prime solutions of linear ternary equations.
\newblock {\em Acta Arith.}, 98(3):293--309, 2001.

\bibitem{Li2003primeandsquare}
Hongze Li.
\newblock The exceptional set for the sum of a prime and a square.
\newblock {\em Acta Math. Hungar.}, 99(1-2):123--141, 2003.

\bibitem{Li2003hybrid}
Hongze Li.
\newblock A hybrid of theorems of {G}oldbach and {P}iatetski-{S}hapiro.
\newblock {\em Acta Arith.}, 107(4):307--326, 2003.

\bibitem{Li2003execptionalset}
Hongze Li.
\newblock A numerical bound for small prime solutions of some binary equations.
\newblock {\em Sci. China Ser. A}, 46(1):48--63, 2003.

\bibitem{LP2010}
Hongze Li and Hao Pan.
\newblock A density version of {V}inogradov's three primes theorem.
\newblock {\em Forum Math.}, 22(4):699--714, 2010.

\bibitem{Li2019}
Huixi Li.
\newblock On the representation of a large integer as the sum of a prime and a
  square-free number with at most three prime divisors.
\newblock {\em Ramanujan J.}, 49(1):141--158, 2019.

\bibitem{Li2023}
Huixi Li.
\newblock Additive representations of natural numbers.
\newblock {\em Ramanujan J.}, 60(4):999--1024, 2023.

\bibitem{LWWY2023}
Huixi {Li}, Biao {Wang}, Chunlin {Wang}, and Shaoyun {Yi}.
\newblock {On the Covering Systems of Polynomial Rings Over Finite Fields}.
\newblock {\em arXiv e-prints}, page arXiv:2308.05378, August 2023.

\bibitem{LWY2023}
Huixi {Li}, Biao {Wang}, and Shaoyun {Yi}.
\newblock {On the minimum modulus problem in number fields}.
\newblock {\em arXiv e-prints}, page arXiv:2302.05946, February 2023.

\bibitem{LZ2022}
Xiaotian Li and Wenguang Zhai.
\newblock The three primes theorem with primes in the intersection of two
  {P}iatetski-{S}hapiro sets.
\newblock {\em Acta Math. Hungar.}, 168(1):228--245, 2022.

\bibitem{Lichtman2022}
Jared~Duker Lichtman.
\newblock Averages of the {M}\"{o}bius function on shifted primes.
\newblock {\em Q. J. Math.}, 73(2):729--757, 2022.

\bibitem{Linnik1951}
Yu.~V. Linnik.
\newblock Prime numbers and powers of two.
\newblock {\em Trudy Mat. Inst. Steklov.}, 38:152--169, 1951.

\bibitem{Linnik1952}
Yu.~V. Linnik.
\newblock Some conditional theorems concerning the binary {G}oldbach problem.
\newblock {\em Izvestiya Akad. Nauk SSSR. Ser. Mat.}, pages 503--520, 1952.

\bibitem{Linnik1953}
Yu.~V. Linnik.
\newblock Addition of prime numbers with powers of one and the same number.
\newblock {\em Mat. Sbornik N.S.}, 32(74):3--60, 1953.

\bibitem{Liu1998}
Jianya Liu.
\newblock The {G}oldbach-{V}inogradov theorem with three primes in a thin
  subset.
\newblock {\em Chinese Ann. Math. Ser. B}, 19(4):479--488, 1998.

\bibitem{LLW1998_1}
Jianya Liu, Mingchit Liu, and Tianze Wang.
\newblock The number of powers of {$2$} in a representation of large even
  integers. {I}.
\newblock {\em Sci. China Ser. A}, 41(4):386--398, 1998.

\bibitem{LLW1998_2}
Jianya Liu, Mingchit Liu, and Tianze Wang.
\newblock The number of powers of {$2$} in a representation of large even
  integers. {II}.
\newblock {\em Sci. China Ser. A}, 41(12):1255--1271, 1998.

\bibitem{LLW1999}
Jianya Liu, Mingchit Liu, and Tianze Wang.
\newblock On the almost {G}oldbach problem of {L}innik.
\newblock {\em J. Th\'{e}or. Nombres Bordeaux}, 11(1):133--147, 1999.

\bibitem{LT2005}
Jianya Liu and Kai-Man Tsang.
\newblock Small prime solutions of ternary linear equations.
\newblock {\em Acta Arith.}, 118(1):79--100, 2005.

\bibitem{Liu1987}
Mingchit Liu.
\newblock An improved bound for prime solutions of some ternary equations.
\newblock {\em Math. Z.}, 194(4):573--583, 1987.

\bibitem{LT1989}
Mingchit Liu and Kaiman Tsang.
\newblock Small prime solutions of linear equations.
\newblock In {\em Th\'{e}orie des nombres ({Q}uebec, {PQ}, 1987)}, pages
  595--624. de Gruyter, Berlin, 1989.

\bibitem{LT1991}
Mingchit Liu and Kaiman Tsang.
\newblock Small prime solutions of some additive equations.
\newblock {\em Monatsh. Math.}, 111(2):147--169, 1991.

\bibitem{LW1998}
Mingchit Liu and Tianze Wang.
\newblock A numerical bound for small prime solutions of some ternary linear
  equations.
\newblock {\em Acta Arith.}, 86(4):343--383, 1998.

\bibitem{LZ1997}
Mingchit Liu and Tao Zhan.
\newblock The {G}oldbach problem with primes in arithmetic progressions.
\newblock In {\em Analytic number theory ({K}yoto, 1996)}, volume 247 of {\em
  London Math. Soc. Lecture Note Ser.}, pages 227--251. Cambridge Univ. Press,
  Cambridge, 1997.

\bibitem{LL2011}
Zhixin Liu and Guangshi L\"{u}.
\newblock Density of two squares of primes and powers of 2.
\newblock {\em Int. J. Number Theory}, 7(5):1317--1329, 2011.

\bibitem{LX2021}
Zhixin Liu and Mengyuan Xue.
\newblock The sum of a prime and a {F}ibonacci number.
\newblock {\em Int. J. Number Theory}, 17(8):1815--1823, 2021.

\bibitem{LY1981}
Shituo Lou and Qi~Yao.
\newblock The exceptional set of {G}oldbach numbers in a short interval.
\newblock {\em Acta Math. Sinica}, 24(2):269--282, 1981.

\bibitem{Lv2007}
Guangshi L\"{u}.
\newblock On {R}omanoff's constant and its generalized problem.
\newblock {\em Adv. Math. (China)}, 36(1):94--100, 2007.

\bibitem{LS2013}
Guangshi L\"{u} and Haiwei Sun.
\newblock The ternary {G}oldbach-{V}inogradov theorem with almost equal primes
  from the {B}eatty sequence.
\newblock {\em Ramanujan J.}, 30(2):153--161, 2013.

\bibitem{Lu2010}
Wenchao Lu.
\newblock Exceptional set of {G}oldbach number.
\newblock {\em J. Number Theory}, 130(10):2359--2392, 2010.

\bibitem{Martin2022}
Kimball Martin.
\newblock Refined {G}oldbach conjectures with primes in progressions.
\newblock {\em Exp. Math.}, 31(1):226--232, 2022.

\bibitem{Matomaki2008}
Kaisa Matom\"{a}ki.
\newblock On the exceptional set in {G}oldbach's problem in short intervals.
\newblock {\em Monatsh. Math.}, 155(2):167--189, 2008.

\bibitem{MMS2017}
Kaisa Matom\"{a}ki, James Maynard, and Xuancheng Shao.
\newblock Vinogradov's theorem with almost equal summands.
\newblock {\em Proc. Lond. Math. Soc. (3)}, 115(2):323--347, 2017.

\bibitem{MRT2019}
Kaisa Matom\"{a}ki, Maksym Radziwi\l~\l, and Terence Tao.
\newblock Correlations of the von {M}angoldt and higher divisor functions {I}.
  {L}ong shift ranges.
\newblock {\em Proc. Lond. Math. Soc. (3)}, 118(2):284--350, 2019.

\bibitem{Meng2007}
Xianmeng Meng.
\newblock A mean value theorem on the binary {G}oldbach problem and its
  application.
\newblock {\em Monatsh. Math.}, 151(4):319--332, 2007.

\bibitem{Meng2009}
Xianmeng Meng.
\newblock On linear equations with prime variables of special type.
\newblock {\em J. Number Theory}, 129(10):2504--2518, 2009.

\bibitem{MW2006}
Xianmeng Meng and Mingqiang Wang.
\newblock A hybrid of theorems of {G}oldbach and {P}iatetski-{S}hapiro.
\newblock {\em Chinese Ann. Math. Ser. B}, 27(3):341--352, 2006.

\bibitem{Mikawa1992}
Hiroshi Mikawa.
\newblock On the exceptional set in {G}oldbach's problem.
\newblock {\em Tsukuba J. Math.}, 16(2):513--543, 1992.

\bibitem{Mikawa1993}
Hiroshi Mikawa.
\newblock On the intervals between consecutive numbers that are sums of two
  primes.
\newblock {\em Tsukuba J. Math.}, 17(2):443--453, 1993.

\bibitem{Mikawa1993_2}
Hiroshi Mikawa.
\newblock On the sum of a prime and a square.
\newblock {\em Tsukuba J. Math.}, 17(2):299--310, 1993.

\bibitem{Mitsui1960}
Takayoshi Mitsui.
\newblock On the {G}oldbach problem in an algebraic number field. {I}, {II}.
\newblock {\em J. Math. Soc. Japan}, 12:290--324, 325--372, 1960.

\bibitem{MV1975}
H.~L. Montgomery and R.~C. Vaughan.
\newblock The exceptional set in {G}oldbach's problem.
\newblock {\em Acta Arith.}, 27:353--370, 1975.

\bibitem{Mozzochi1980}
C.~J. Mozzochi.
\newblock An analytic sufficiency condition for {G}oldbach's conjecture with
  minimal redundancy.
\newblock {\em Kyungpook Math. J.}, 20(1):1--9, 1980.

\bibitem{Mozzochi1981}
C.~J. Mozzochi.
\newblock An analytic sufficiency condition for {G}oldbach's conjecture with
  minimal redundancy. {II}.
\newblock {\em Kyungpook Math. J.}, 21(1):5--8, 1981.

\bibitem{Mozzochi2014}
C.~J. Mozzochi.
\newblock A comparison of sufficiency condtions for the goldbach and the twin
  primes conjectures.
\newblock {\em Advances in Pure Mathematics}, 4:157--170, 2014.

\bibitem{MB1978}
C.~J. Mozzochi and R.~Balasubramanian.
\newblock Some comments on {G}oldbach's conjecture.
\newblock {\em Report No. 11, Mittag-Leffler Institute}, 1978.

\bibitem{Na1986}
Jisheng Na.
\newblock An estimation related to {H}ua's attempt at {G}oldbach problem.
\newblock {\em Kexue Tongbao (English Ed.)}, 31(22):1513--1520, 1986.

\bibitem{Nayebi2011}
A.~Nayebi.
\newblock Upper bounds on the solutions to {$n=p+m^2$}.
\newblock {\em Bull. Iranian Math. Soc.}, 37(4):95--108, 2011.

\bibitem{Nielsen2009}
Pace~P. Nielsen.
\newblock A covering system whose smallest modulus is 40.
\newblock {\em J. Number Theory}, 129(3):640--666, 2009.

\bibitem{Encyclopedia}
Encyclopedia of~Mathematics.
\newblock Goldbach-{W}aring problem.
\newblock
  \url{http://encyclopediaofmath.org/index.php?title=Goldbach-Waring_problem&oldid=44620}.
\newblock Accessed: 2023.

\bibitem{OHP2014}
Tom\'{a}s Oliveira~e Silva, Siegfried Herzog, and Silvio Pardi.
\newblock Empirical verification of the even {G}oldbach conjecture and
  computation of prime gaps up to {$4\cdot 10^{18}$}.
\newblock {\em Math. Comp.}, 83(288):2033--2060, 2014.

\bibitem{Pan1959}
Chengdong Pan.
\newblock Some new results in the additive theory of pirme numbers.
\newblock {\em Acta Math. Sinica}, 9:315--329, 1959.

\bibitem{Pan1982}
Chengdong Pan.
\newblock A new attempt on {G}oldbach conjecture.
\newblock {\em Chinese Ann. Math.}, 3(4):555--560, 1982.

\bibitem{PDW1975}
Chengdong Pan, Xiaqi Ding, and Yuan Wang.
\newblock On the representation of every large even integer as a sum of a prime
  and an almost prime.
\newblock {\em Sci. Sinica}, 18(5):599--610, 1975.

\bibitem{PP1989}
Chengdong Pan and Chengbiao Pan.
\newblock On estimations of trigonometric sums over primes in short intervals.
  {II}.
\newblock {\em Sci. China Ser. A}, 32(6):641--653, 1989.

\bibitem{PP1990}
Chengdong Pan and Chengbiao Pan.
\newblock On estimations of trigonometric sums over primes in short intervals.
  {III}.
\newblock {\em Chinese Ann. Math. Ser. B}, 11(2):138--147, 1990.

\bibitem{Paran2020}
Elad Paran.
\newblock Twin-prime and {G}oldbach theorems for {$\mathbb{Z}[[x]]$}.
\newblock {\em J. Number Theory}, 213:453--461, 2020.

\bibitem{Peneva2000}
T.~P. Peneva.
\newblock On the ternary {G}oldbach problem with primes {$p_i$} such that
  {$p_i+2$} are almost-primes.
\newblock {\em Acta Math. Hungar.}, 86(4):305--318, 2000.

\bibitem{Peneva2001}
T.~P. Peneva.
\newblock On the exceptional set for {G}oldbach's problem in short intervals.
\newblock {\em Monatsh. Math.}, 132(1):49--65, 2001.

\bibitem{PP1993}
A.~Perelli and J.~Pintz.
\newblock On the exceptional set for {G}oldbach's problem in short intervals.
\newblock {\em J. London Math. Soc. (2)}, 47(1):41--49, 1993.

\bibitem{PP1995}
A.~Perelli and J.~Pintz.
\newblock Hardy-{L}ittlewood numbers in short intervals.
\newblock {\em J. Number Theory}, 54(2):297--308, 1995.

\bibitem{PZ1995}
A.~Perelli and A.~Zaccagnini.
\newblock On the sum of a prime and a {$k$}th power.
\newblock {\em Izv. Ross. Akad. Nauk Ser. Mat.}, 59(1):185--200, 1995.

\bibitem{Pintz2006}
J.~Pintz.
\newblock A note on {R}omanov's constant.
\newblock {\em Acta Math. Hungar.}, 112(1-2):1--14, 2006.

\bibitem{Pintz2012}
J.~Pintz.
\newblock The bounded gap conjecture and bounds between consecutive {G}oldbach
  numbers.
\newblock {\em Acta Arith.}, 155(4):397--405, 2012.

\bibitem{Pintz2018_1}
J.~{Pintz}.
\newblock {A new explicit formula in the additive theory of primes with
  applications I. The explicit formula for the Goldbach and Generalized Twin
  Prime Problems}.
\newblock {\em arXiv e-prints}, page arXiv:1804.05561, April 2018.

\bibitem{Pintz2018_2}
J.~{Pintz}.
\newblock {A new explicit formula in the additive theory of primes with
  applications II. The exceptional set in Goldbach's problem}.
\newblock {\em arXiv e-prints}, page arXiv:1804.09084, April 2018.

\bibitem{Pintz2023}
J.~Pintz.
\newblock {A new explicit formula in the additive theory of primes with
  applications I. The explicit formula for the Goldbach and Generalized Twin
  Prime Problems}.
\newblock {\em Acta Arith.}, pages 10.4064/aa220728--31--3, 2023.

\bibitem{PR2003}
J.~Pintz and I.~Z. Ruzsa.
\newblock On {L}innik's approximation to {G}oldbach's problem. {I}.
\newblock {\em Acta Arith.}, 109(2):169--194, 2003.

\bibitem{PR2020}
J.~Pintz and I.~Z. Ruzsa.
\newblock On {L}innik's approximation to {G}oldbach's problem. {II}.
\newblock {\em Acta Math. Hungar.}, 161(2):569--582, 2020.

\bibitem{Pollack2011}
Paul Pollack.
\newblock On polynomial rings with a {G}oldbach property.
\newblock {\em Amer. Math. Monthly}, 118(1):71--77, 2011.

\bibitem{Polyakov1990}
I.~V. Polyakov.
\newblock Addition of a prime number and the square of an integer.
\newblock {\em Mat. Zametki}, 47(4):90--99, 1990.

\bibitem{Polymath2014}
D.~H.~J. Polymath.
\newblock Variants of the {S}elberg sieve, and bounded intervals containing
  many primes.
\newblock {\em Res. Math. Sci.}, 1:Art. 12, 83, 2014.

\bibitem{PW2022}
Carl Pomerance and Andreas Weingartner.
\newblock On primes and practical numbers.
\newblock {\em Ramanujan J.}, 57(3):981--1000, 2022.

\bibitem{Qin1996}
Hourong Qin.
\newblock Fermat's problem and {G}oldbach's problem over {$M_n{\bf Z}$}.
\newblock {\em Linear Algebra Appl.}, 236:131--135, 1996.

\bibitem{RV2017}
M.~Ram~Murty and Akshaa Vatwani.
\newblock Twin primes and the parity problem.
\newblock {\em J. Number Theory}, 180:643--659, 2017.

\bibitem{Ramachandra1973}
K.~Ramachandra.
\newblock On the number of {G}oldbach numbers in small intervals.
\newblock {\em J. Indian Math. Soc. (N.S.)}, 37:157--170, 1973.

\bibitem{Richstein2001}
J\"{o}rg Richstein.
\newblock Verifying the {G}oldbach conjecture up to {$4\cdot 10^{14}$}.
\newblock {\em Math. Comp.}, 70(236):1745--1749, 2001.

\bibitem{Rivat1992}
J.~Rivat.
\newblock Autour d'un th$\acute{e}$orem$\grave{e}$ de {P}iatetski-{S}hapiro.
\newblock {\em Thesis}, Paris XI, Orsay, 1992.

\bibitem{Romani1983}
F.~Romani.
\newblock Computations concerning primes and powers of two.
\newblock {\em Calcolo}, 20(3):319--336 (1984), 1983.

\bibitem{Romanoff1934}
N.~P. Romanoff.
\newblock \"{U}ber einige {S}\"{a}tze der additiven {Z}ahlentheorie.
\newblock {\em Math. Ann.}, 109(1):668--678, 1934.

\bibitem{Ross1975}
P.~M. Ross.
\newblock On {C}hen's theorem that each large even number has the form
  {$p_{1}+p_{2}$} or {$p_{1}+p_{2}p_{3}$}.
\newblock {\em J. London Math. Soc. (2)}, 10(4):500--506, 1975.

\bibitem{Ross1978}
P.~M. Ross.
\newblock A short intervals result in additive prime number theory.
\newblock {\em J. London Math. Soc. (2)}, 17(2):219--227, 1978.

\bibitem{Saidak2006}
Filip Saidak.
\newblock On {G}oldbach's conjecture for integer polynomials.
\newblock {\em Amer. Math. Monthly}, 113(6):541--545, 2006.

\bibitem{SV1993}
Saverio Salerno and Antonio Vitolo.
\newblock {$p+2=P_2$} in short intervals.
\newblock {\em Note Mat.}, 13(2):309--328, 1993.

\bibitem{SC2021}
H.~M. Saliba and V.~N. Chubarikov.
\newblock On some additive problems of {G}oldbach's type.
\newblock {\em Chebyshevski\u{\i} Sb.}, 22(3):179--195, 2021.

\bibitem{Salmensuu2022}
Juho Salmensuu.
\newblock The {G}oldbach conjecture with summands in arithmetic progressions.
\newblock {\em Q. J. Math.}, 73(4):1375--1401, 2022.

\bibitem{Shao2014}
Xuancheng Shao.
\newblock A density version of the {V}inogradov three primes theorem.
\newblock {\em Duke Math. J.}, 163(3):489--512, 2014.

\bibitem{Shen2016}
Quanli Shen.
\newblock The ternary {G}oldbach problem with primes in positive density sets.
\newblock {\em J. Number Theory}, 168:334--345, 2016.

\bibitem{SW2010}
Igor~E. Shparlinski and Arne Winterhof.
\newblock Partitions into two {L}ehmer numbers.
\newblock {\em Monatsh. Math.}, 160(4):429--441, 2010.

\bibitem{Stein1856}
Moritz~A. Stein.
\newblock Sur un assertion de {G}oldbach relative aux nombres impairs.
\newblock {\em Nouvelles Annales Math.}, 15:23--24, 1856.

\bibitem{Sun2000}
Zhiwei Sun.
\newblock On integers not of the form {$\pm p^a\pm q^b$}.
\newblock {\em Proc. Amer. Math. Soc.}, 128(4):997--1002, 2000.

\bibitem{Suzuki2017}
Yuta Suzuki.
\newblock A remark on the conditional estimate for the sum of a prime and a
  square.
\newblock {\em Funct. Approx. Comment. Math.}, 57(1):61--76, 2017.

\bibitem{Tolev1999}
D.~I. Tolev.
\newblock Arithmetic progressions of prime-almost-prime twins.
\newblock {\em Acta Arith.}, 88(1):67--98, 1999.

\bibitem{Tolev2000_1}
D.~I. Tolev.
\newblock Additive problems with prime numbers of special type.
\newblock {\em Acta Arith.}, 96(1):53--88, 2000.

\bibitem{Tolev2000_2}
D.~I. Tolev.
\newblock Representations of large integers as sums of two primes of special
  type.
\newblock In {\em Algebraic number theory and {D}iophantine analysis ({G}raz,
  1998)}, pages 485--495. de Gruyter, Berlin, 2000.

\bibitem{DT2023}
Ognian {Trifonov} and Jack {Dalton}.
\newblock {Representing positive integers as a sum of a squarefree number and a
  small prime}.
\newblock {\em arXiv e-prints}, page arXiv:2301.12585, January 2023.

\bibitem{Vaserstein1989}
L.~N. Vaserstein.
\newblock Noncommutative number theory.
\newblock In {\em Algebraic {$K$}-theory and algebraic number theory
  ({H}onolulu, {HI}, 1987)}, volume~83 of {\em Contemp. Math.}, pages 445--449.
  Amer. Math. Soc., Providence, RI, 1989.

\bibitem{Vaughan1972}
R.~C. Vaughan.
\newblock On {G}oldbach's problem.
\newblock {\em Acta Arith.}, 22:21--48, 1972.

\bibitem{Vaughan2014}
R.~C. Vaughan.
\newblock The general {G}oldbach problem with {B}eatty primes.
\newblock {\em Ramanujan J.}, 34(3):347--359, 2014.

\bibitem{Vaughan2016}
R.~C. Vaughan.
\newblock Goldbach's conjectures: a historical perspective.
\newblock In {\em Open problems in mathematics}, pages 479--520. Springer,
  [Cham], 2016.

\bibitem{Vinogradov1985}
A.~I. Vinogradov.
\newblock The binary {H}ardy-{L}ittlewood problem.
\newblock {\em Acta Arith.}, 46(1):33--56, 1985.

\bibitem{Wang1992}
Jun Wang.
\newblock Goldbach's problem in the ring {$M_n({\bf Z})$}.
\newblock {\em Amer. Math. Monthly}, 99(9):856--857, 1992.

\bibitem{WC2023}
Ruijing Wang and Yonggao Chen.
\newblock On the sum of a {F}ibonacci number and a prime.
\newblock {\em Int. J. Number Theory}, 19(4):873--889, 2023.

\bibitem{Wang1995}
Tianze Wang.
\newblock On the exceptional set for the equation {$n=p+k^2$}.
\newblock {\em Acta Math. Sinica (N.S.)}, 11(2):156--167, 1995.

\bibitem{Wang1999}
Tianze Wang.
\newblock On {L}innik's almost {G}oldbach theorem.
\newblock {\em Sci. China Ser. A}, 42(11):1155--1172, 1999.

\bibitem{Wang1982}
Yuan Wang.
\newblock The {G}oldbach conjecture.
\newblock {\em Math. Medley}, 10(1):1--4, 1982.

\bibitem{Wong1996}
K.C. Wong.
\newblock Contributions to analytic number theory.
\newblock {\em Ph.D. Thesis}, Cardiff, 1996.

\bibitem{Wu1993}
Jie Wu.
\newblock Th\'{e}or\`emes g\'{e}n\'{e}ralis\'{e}s de {B}ombieri-{V}inogradov
  dans les petits intervalles.
\newblock {\em Quart. J. Math. Oxford Ser. (2)}, 44(173):109--128, 1993.

\bibitem{Wu1994}
Jie Wu.
\newblock Sur l'\'{e}quation {$p+2=P_2$} dans les petits intervalles.
\newblock {\em J. London Math. Soc. (2)}, 49(1):61--72, 1994.

\bibitem{Wu2004}
Jie Wu.
\newblock Chen's double sieve, {G}oldbach's conjecture and the twin prime
  problem.
\newblock {\em Acta Arith.}, 114(3):215--273, 2004.

\bibitem{Wu2008}
Jie Wu.
\newblock Chen's double sieve, {G}oldbach's conjecture and the twin prime
  problem. {II}.
\newblock {\em Acta Arith.}, 131(4):367--387, 2008.

\bibitem{Xylouris2011}
Triantafyllos Xylouris.
\newblock On the least prime in an arithmetic progression and estimates for the
  zeros of {D}irichlet {$L$}-functions.
\newblock {\em Acta Arith.}, 150(1):65--91, 2011.

\bibitem{Yamada2015}
Tomohiro {Yamada}.
\newblock {Explicit Chen's theorem}.
\newblock {\em arXiv e-prints}, page arXiv:1511.03409, November 2015.

\bibitem{Yao1982}
Qi~Yao.
\newblock The exceptional set of {G}oldbach numbers in a short interval.
\newblock {\em Acta Math. Sinica}, 25(3):315--322, 1982.

\bibitem{Yau2021}
K.~Yau.
\newblock Representation of an integer as the sum of a prime in arithmetic
  progression and a square-free integer.
\newblock {\em Funct. Approx. Comment. Math.}, 64(1):77--108, 2021.

\bibitem{Zhan1991}
Tao Zhan.
\newblock On the representation of large odd integer as a sum of three almost
  equal primes.
\newblock {\em Acta Math. Sinica (N.S.)}, 7(3):259--272, 1991.

\bibitem{Zhan1995}
Tao Zhan.
\newblock A generalization of the {G}oldbach-{V}inogradov theorem.
\newblock {\em Acta Arith.}, 71(2):95--106, 1995.

\bibitem{Zhang2010}
Shaohua Zhang.
\newblock Goldbach conjecture and the least prime number in an arithmetic
  progression.
\newblock {\em C. R. Math. Acad. Sci. Paris}, 348(5-6):241--242, 2010.

\bibitem{Zhang2014}
Yitang Zhang.
\newblock Bounded gaps between primes.
\newblock {\em Ann. of Math. (2)}, 179(3):1121--1174, 2014.

\bibitem{Zhang2022}
Yitang {Zhang}.
\newblock {Discrete mean estimates and the Landau-Siegel zero}.
\newblock {\em arXiv e-prints}, page arXiv:2211.02515, November 2022.

\end{thebibliography}

\end{document}